\documentclass[12pt]{amsart}
\usepackage[normalem]{ulem} % provides \sout for strikethrough

\usepackage[utf8]{inputenc}

\usepackage{graphicx}
\usepackage{verbatim}
\usepackage{textcomp}
\usepackage{amssymb}
\usepackage{cite}
\usepackage{color}
\usepackage[all]{xy}
\usepackage{graphicx}
\usepackage{color}
\usepackage{hyperref}
\hypersetup{
    colorlinks,
    citecolor=black,
    filecolor=black,
    linkcolor=black,
    urlcolor=black
}
\usepackage{tikz}
\usepackage[all]{xy}
\usepackage{graphicx}
\usetikzlibrary{backgrounds}
\newtheorem{theorem}{Theorem}

\allowdisplaybreaks

\theoremstyle{definition}
\newtheorem{definition}[theorem]{Definition}

\theoremstyle{remark}
\newtheorem{remark}[theorem]{Remark}

\numberwithin{equation}{section}

\newtheorem*{theorem*}{Theorem}
{{\sc Proof of Lemma~\ref{tri1}.}}%
{{\qed} \\}
{{\sc Proof of Theorem~\ref{regularity}.}}%
{{\qed} \\}
{{\sc Proof of Theorem~\ref{main}.}}%
{{\qed} \\}
\newenvironment{proofof(i)}%
    {{\sc Proof of $(i)$.}}%
  {{\qed} \\}  
  
  \newenvironment{proofof(iv)}%
    {{\sc Proof of $(iv)$.}}%
  {{\qed} \\}  
\newcommand{\R}{\mathbb R}

\newcommand{\C}{\mathbb C}

\title{Harmonic Maps in Singular Geometry and Rigidity}

\thanks{}

\author[Daskalopoulos]{Georgios Daskalopoulos}
\address{Brown University\\
Department of Mathematics\\
Providence, RI}%  21218}
\email{daskal@math.brown.edu}

\author[Mese]{Chikako Mese}
\address{Johns Hopkins University\\
Department of Mathematics\\
Baltimore, MD}%  21218}
\email{mese@math.jhu.edu}

\begin{document}
\maketitle

\begin{abstract}
This survey reviews results on harmonic maps into spaces of non-positive curvature, with a focus on targets that lack smooth structure. More precisely, we consider targets that are complete metric spaces with non-positive curvature in the sense of Alexandrov, commonly referred to as NPC (non-positively curved) or CAT(0) spaces. We discuss applications of harmonic maps to rigidity phenomena, including generalizations of Margulis superrigidity and the holomorphic rigidity of Teichm\"uller space.
Our approach relies heavily on the regularity theory of harmonic maps to non-smooth targets, enabling differential-geometric techniques to be employed in the absence of any smooth structure on the target.

\end{abstract}
\tableofcontents  % This generates the table of contents

\section{Background on harmonic maps}
\subsection{Smooth targets}\label{smoothharma}

 Given a smooth map \( u: M \to N \) between Riemannian manifolds \((M, g)\) and \((N, h)\), the energy of \( u \) is defined by
\[
E(u) = \frac{1}{2} \int_M |du|^2 \, d\mathrm{vol}_g,
\]
where
\[
|du|^2 = \sum_{i=1}^m h(du(e_i), du(e_i)),
\]
and \( \{e_1, \dots, e_m\} \) is any local orthonormal frame for the tangent bundle \( TM \). The scalar quantity \( |du|^2 \), called the \emph{energy density}, is independent of the choice of orthonormal basis.

Critical points of the energy functional are called \emph{harmonic maps}. These generalize both harmonic functions (when \( \dim N = 1 \)) and geodesics (when \( \dim M = 1 \)). The Euler-Lagrange equation corresponding to the energy functional gives rise to a second-order semilinear elliptic partial differential equation known as the \emph{harmonic map equation}:
\begin{equation} \label{hme}
\Delta u^i + \sum_{\alpha, \beta} \sum_{j,k} g^{\alpha \beta} \Gamma^i_{jk} \circ u \frac{\partial u^j}{\partial x^\alpha} \frac{\partial u^k}{\partial x^\beta} = 0, \quad i = 1, \dots, m
\end{equation}
where \( \Gamma_{ij}^k \)  denote the Christoffel symbols of \( N \).

One of the central questions is the study of existence and regularity of harmonic maps. In the 1960's, foundational work by Eells-Sampson \cite{eells-sampson} showed that if \( N \) has non-positive sectional curvature and \( u_0: M \to N \) is of finite energy, then there exists a smooth harmonic map $u$ in the  homotopy class of $u_0$. The map $u$ is obtained from the initial map $u_0$ by the \emph{harmonic map heat flow}, a parabolic flow that decreases energy and converges (under appropriate curvature conditions) to a harmonic map. 

Many of the important  applications of harmonic maps to geometry and topology are well known and beyond the scope of this note to discuss. Here we will concentrate in its role to the rigidity theory of discrete groups and locally symmetric spaces. For instance, in the context of Mostow rigidity \cite{mostow},  uniqueness of harmonic maps \cite{hartman} implies that the harmonic map  must in fact be an isometry. This  illustrates a general strategy: one seeks a harmonic map and then uses the geometric properties of the domain and/or target to deduce that the map satisfies stronger conditions such as being totally geodesic, an isometry, or even constant.  This method has been particularly successful in yielding an alternative proof and various extensions of the Margulis superrigidity theorem (cf. \cite{margulis}), although a harmonic map proof of Mostow's theorem has yet to be found.

In order to illustrate the method, we will briefly describe the Eells-Sampson Bochner formula:
\begin{theorem}\label{BochnerES}
Let $u: M \to N$ be a harmonic map between compact Riemannian manifolds, and $(e_\alpha)$ an orthonormal frame for $TM$. Then 
\[
\Delta |du|^2 
= |\nabla d u|^2 
+ \tfrac{1}{2} \ \langle d u(Ric^M(e_\alpha)), d u(e_\alpha) \rangle 
- \tfrac{1}{2}  \langle R^N(d u(e_\alpha), d u(e_\beta))d u(e_\beta), d u(e_\alpha) \rangle.
\]
\end{theorem}

In the above, $\nabla d u$ denotes the Hessian of $u$ while $Ric^M$ and $R^N$ are the Ricci and Riemannian tensors on the domain and target, respectively.
Assume now that $Ric^M \geq 0$,  and $R^N \leq 0$. Integrating over $M$, and noting that the integral of the Laplacian is zero, we conclude that $\nabla d u=0$; in other words, $u$  is totally geodesic. This provides the first and simplest illustration of how, under suitable curvature conditions, a harmonic map exhibits rigidity.

 In \cite{siu}, Siu  established strong rigidity  for certain classes of compact  K\"ahler manifolds. Specifically, he proved that any smooth map between compact K\"ahler manifolds  is homotopic to a  holomorphic or anti-holomorphic map, provided the target has suitably strong negative curvature. The key is again the discovery of a  Bochner formula, which implies first that the harmonic map is pluriharmonic (i.e.~its restriction to any complex disc is also harmonic). Further analysis shows that the map is, in fact, holomorphic or antiholomorphic. Siu's approach paved the way for subsequent rigidity results via harmonic maps, obtained by Corlette \cite{corlette}, Jost-Yau \cite{jost-yau}, Mok-Siu-Yeung \cite{mok-siu-yeung}, and others.

\subsection{NPC space targets}
\label{npctarget}

In \cite{korevaar-schoen} and \cite{korevaar-schoen2}, Korevaar and Schoen consider maps \( u: \Omega \to X \) from an \( m \)-dimensional Riemannian domain (i.e., a connected, open, subset of a Riemannian manifold \( (M,g) \), with compact closure in \( M \)) to an NPC space.

\begin{figure}[h]
\begin{center}
\includegraphics[width=0.6\textwidth]{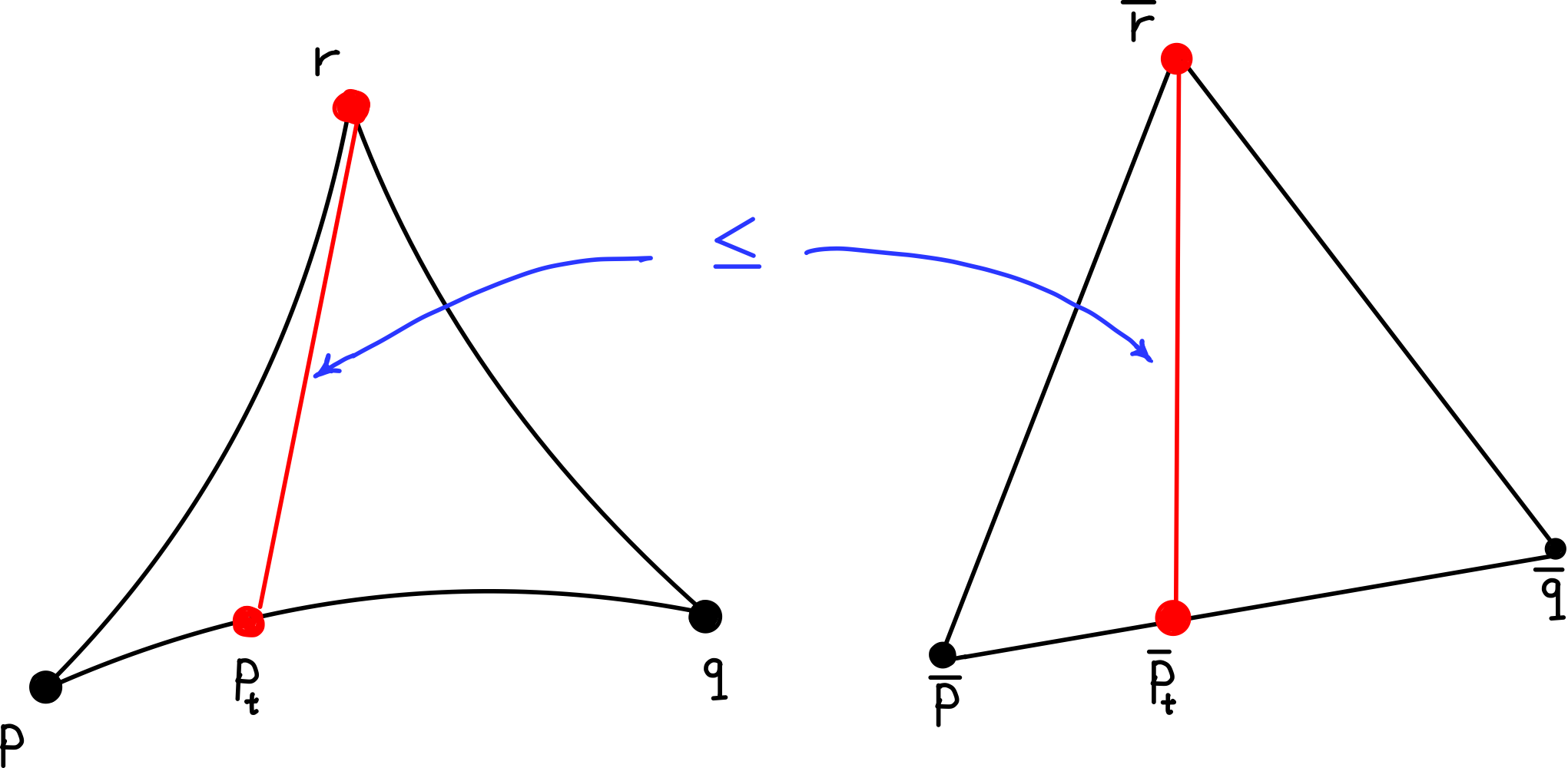}
\caption{Triangle comparison}%date unknown
\label{Triangle}
\end{center}
\end{figure}

\begin{definition}
A simply connected, complete metric space $(X,d)$ is said to be an NPC space if the following conditions hold: 

\begin{itemize}
\item[(1)] It is a  geodesic space; i.e.~for every pair of points $x_0, x_1 \in X$, there exists a continuous map $\gamma:[0,1] \to X$ such that  

\begin{itemize}
\item[(i)] $\gamma(0)=x_0$ and $\gamma(1)=x_1$
\item[(ii)] For all $t, s \in [0,1]$, $d(\gamma(t),\gamma(s)) =|t-s|d(x_0,x_1)$.
\end{itemize} 

\item[(2)] It  has nonpositive curvature in the sense of triangle comparison; i.e.~a geodesic triangle in $X$ with vertices $p,q,r$ is no thicker than a comparison geodesic triangle  in Euclidean plane ${\mathbb E}^2$.  More precisely, \emph{if  $\bar{p}, \bar{q}, \bar{r} \in {\mathbb{E}}^2$ are such that }
 \[
d(p,q)=|\bar{p}-\bar{q}|, \,\, d(q,r)=|\bar{q}-\bar{r}|, \,\, d(r,p)=|\bar{r}-\bar{p}|,
\]
{\em then }
\[
d(p_t,r)\le |\big((1-t)\bar{p}+t\bar{q}\big) - \bar{r}|
\]
\emph{where $t \mapsto p_t$ for $ t\in [0,1]$  is the constant speed parameterization of the  geodesic  from $p$ to $q$.}
See Figure~\ref{Triangle}.  
\end{itemize}
\end{definition}
  The study of NPC spaces (more generally, spaces with curvature bounded from above by $\kappa$ known as CAT($\kappa$) spaces) was initiated by  A.~D.~Alexandrov  and further brought to prominence by  Gromov.  

\begin{figure}[h]
\begin{center}
\includegraphics[width=0.2\textwidth]{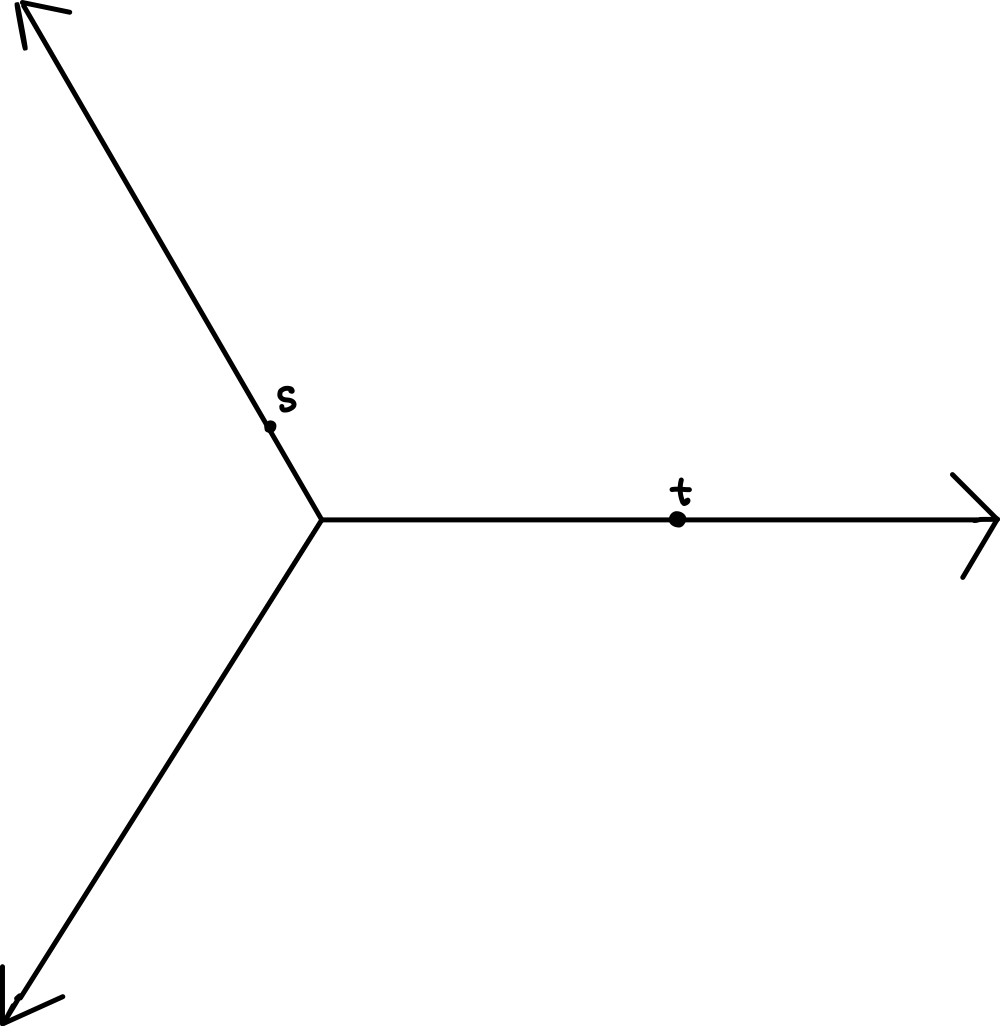}
\caption{tripod}%date unknown
\label{tripod}
\end{center}
\end{figure}

Hadamard manifolds, i.e.~complete, simply connected Riemannian manifolds with sectional curvature bounded above by zero, are  examples of NPC spaces. A simple example of an NPC space that is not a manifold is the tripod (Figure \ref{tripod}), obtained by gluing three copies of the interval $[0,\infty)$ at their origins, called the vertex and denoted by $V$. The distance between two points $t$ and $s$ lying in different copies of $[0,\infty)$ is defined to be $t+s$. Observe that any pair of points in a tripod lies in a totally geodesic, isometrically embedded copy of $\mathbb{R}$ in $T$. More generally, one can construct a $k$-pod by gluing $k$ copies of $[0,\infty)$ in the same way.

A higher dimensional generalization of trees is an \emph{$n$-dimensional Euclidean building} $X$. Like the $k$-pod, any two points $p, q \in X$ lie in the image of a totally geodesic, isometric embedding of an $n$-dimensional Euclidean space called \emph{apartment}. Euclidean buildings are built from apartments  glued together along convex polyhedral subsets according to combinatorial rules determined by an affine Weyl  group. In this paper, we adopt the geometric definition of a Euclidean building as given in \cite{kleiner-leeb}.

Other important examples of NPC spaces are \emph{hyperbolic buildings} arising in the study of hyperbolic groups, and the \emph{Weil-Petersson completion} of Teichm\"uller space.

We now briefly describe  the basics of harmonic map theory to NPC spaces, according to \cite{korevaar-schoen} and \cite{korevaar-schoen2}.
To first define \(\varepsilon\)-approximate energy density at $x$, consider the square of the difference quotient at a point \( x \in \Omega \), and integrate it over unit vectors \( V \in S_1 \subset T_xM \) 
\[
e_\varepsilon(x) = \int_{y \in \partial B_{\epsilon}(x)}  \frac{d^2(u(x), u(y))}{\varepsilon^2} 
\frac{d\sigma_{x,\varepsilon}(y)}{\varepsilon^{m-1}}
\]
where $d\sigma_{x,\varepsilon}(y)$ is the $(m-1)$-dimensional surface measure on $\partial B_\epsilon(x)$.
Let \( u \in L^2(\Omega,X) \); that is, for any \( P \in X \),
\[
\int_\Omega d^2(u(x),P)\, d\mathrm{vol}_g(x) < \infty.
\]
For \( \varepsilon > 0 \) and \( f \in C_c(\Omega) \), define
\[
E_\varepsilon^u(f) = \int_\Omega f(x)\, e_\varepsilon(x)\, d\mathrm{vol}_g(x).
\]
We say that \( u \) has \emph{finite energy} (or \( u \in W^{1,2}(\Omega,X) \)) if
\[
E^u := \sup_{\substack{f \in C_c(\Omega)\\0 \leq f \leq 1}} \limsup_{\varepsilon \to 0} E_\varepsilon^u(f) < \infty.
\]

If \( E^u < \infty \), then the measures \( e_\varepsilon(x)\, d\mathrm{vol}_g(x) \) converge weakly (as \( \varepsilon \to 0 \)) to a measure of the form \( e(x)\, d\mathrm{vol}_g(x) \), where \( |\nabla u|^2 := e(x) \) is an integrable function. The \emph{energy} is thus defined by
\[
E^u= \int_\Omega |\nabla u|^2\, d\mathrm{vol}_g.
\]
The map \( u \) is said to be \emph{harmonic} if it is locally energy minimizing; i.e., for every \( x \in \Omega \), there exists \( r > 0 \) such that \( u|_{B_r(x)} \) minimizes energy among all maps \( v: B_r(x) \to X \) with the same trace (i.e.~boundary values) as \( u|_{B_r(x)} \).

Existence and regularity of harmonic maps into an NPC space is established in \cite{gromov-schoen}, \cite{korevaar-schoen, korevaar-schoen2, korevaar-schoen3}.  For example, for the Dirichlet problem: 

\begin{theorem}[\cite{korevaar-schoen} Theorem 2.2 and 2.4.6]
\label{2.4.6}
Let $(\Omega,g)$ be a Lipschitz Riemannian domain and let $(X,d)$ be an
NPC metric space. Let $ \phi \in W^{1,2}(\Omega,X)$ and define
\[
W^{1,2}_{\phi}(\Omega,X)=\{u \in W^{1,2}(\Omega,X):  \mathrm{tr}(u) = \mathrm{tr}(\phi)\}.
\]
Then there exists a unique $u \in W^{1,2}_{\phi}(\Omega,X)$ such 
\[
E^u \leq \inf_{v \in W^{1,2}_\phi(\Omega,X)} E^v.
\]
Moreover, $u$ is a locally Lipschitz continuous function in the interior of $\Omega$ where the local Lipschitz constant is bounded above by a constant dependent on $E^u$, the distance to the boundary $\partial \Omega$ and the  metric $g$.
\end{theorem}

For applications to rigidity problems, it is natural to consider a generalization 
of the homotopy problem, namely the \emph{equivariant} problem.  
Let $M$ be a  manifold with universal cover $\widetilde M$, $X$ a metric space, 
and   
\[
\rho : \pi_1(M) \to \mathsf{Isom}(X)
\]  
 a homomorphism intertwining the  action  by deck transformations on  
$\widetilde M$ with the action by isometries on $X$.  
A map $f : \widetilde M \to X$ is called \emph{$\rho$-equivariant} if  
\[
f(\gamma p) = \rho(\gamma) f(p) \qquad 
\text{for all } \gamma \in \pi_1(M), \ p \in \widetilde M.
\]  
If $f$ is $\rho$-equivariant, the energy density $|df|^2$ on $\widetilde M$ is 
$\pi_1(M)$-invariant and thus descends to a function on $M$. 
We define the \emph{energy} of $f$ by  
\[
E^f:= \int_M |df|^2 \, d\mathrm{vol}_g.
\]  
\begin{definition}
A $\rho$-equivariant map $ u:\widetilde{M} \to X$ is harmonic if $E^u \leq E^f$ for any $\rho$-equivariant map $f: \widetilde M \to X$.
\end{definition}

 \begin{theorem}\label{existence}[Existence, \cite{gromov-schoen} and \cite{korevaar-schoen2}]  
Let $M$ be a finite volume Riemannian manifold, $\Gamma =\pi_1(M)$ and $X$ a locally compact NPC space.  Assume $\rho:\Gamma \to \mathsf{Isom}(X)$  does not fix a point at infinity.  If  there exists a finite energy $\rho$-equivariant map, then there exists 
 a finite energy $\rho$-equivariant harmonic map $u: \widetilde M \rightarrow X$.
 \end{theorem}

\begin{remark} The assumption that $\rho$ does not fix a point at infinity is to  keep an energy minimizing sequence of maps from escaping to the boundary at infinity. 
\end{remark}

Harmonic maps into NPC spaces also enjoy uniqueness properties.  Indeed, by the quadrilateral comparison inequality (a special case  of Reshetnyak's theorem on the convexity of the distance function in CAT($\kappa$) spaces \cite{reshetnyak}), 
\[
d^2(p_t,q_t) \leq (1-t) d^2(p,q)+ t d^2(r,s)-t(1-t) (d(p,s)-d(q,r))^2
\]
for the ordered sequence $\{p,q,r,s\} \subset X$ where $p_t$ (resp.~$q_t$)  is the point on  the unique  geodesic from $p$ to $s$ (resp.~$q$ to $r$) at a distance $td(p,s)$ from $p$ (resp.~$td(q,r)$ from $q$).
This  in turn implies    the convexity of the energy functional
\[
\int_M |\nabla u_t|^2 d\mathrm{vol}_g \leq (1-t) \int_M |\nabla u_0|^2 d\mathrm{vol}_g + t \int_M |\nabla u_1|^2 d\mathrm{vol}_g
 - t(1-t) \int_M |\nabla d(u_0,u_1)|^2 d\mathrm{vol}_g
\]
where  $u_0,u_1 \in W^{1,2}(M,X)$ and $u_t$ is the interpolation map; i.e.~$u_t(x)$ equals the point on the unique geodesic between  $u_0(x)$ and $u_1(x)$ at a distance $td(u_0(x),u_1(x))$ from $u_0(x)$ for $t \in [0,1]$.  Thus, if $u_0$ and $u_1$ are energy minimizing, then the map $u_t$ is also energy minimizing and $d(u_0,u_1)$ is a constant.

\begin{figure}[h]
\begin{center}
\includegraphics[width=0.4\textwidth]{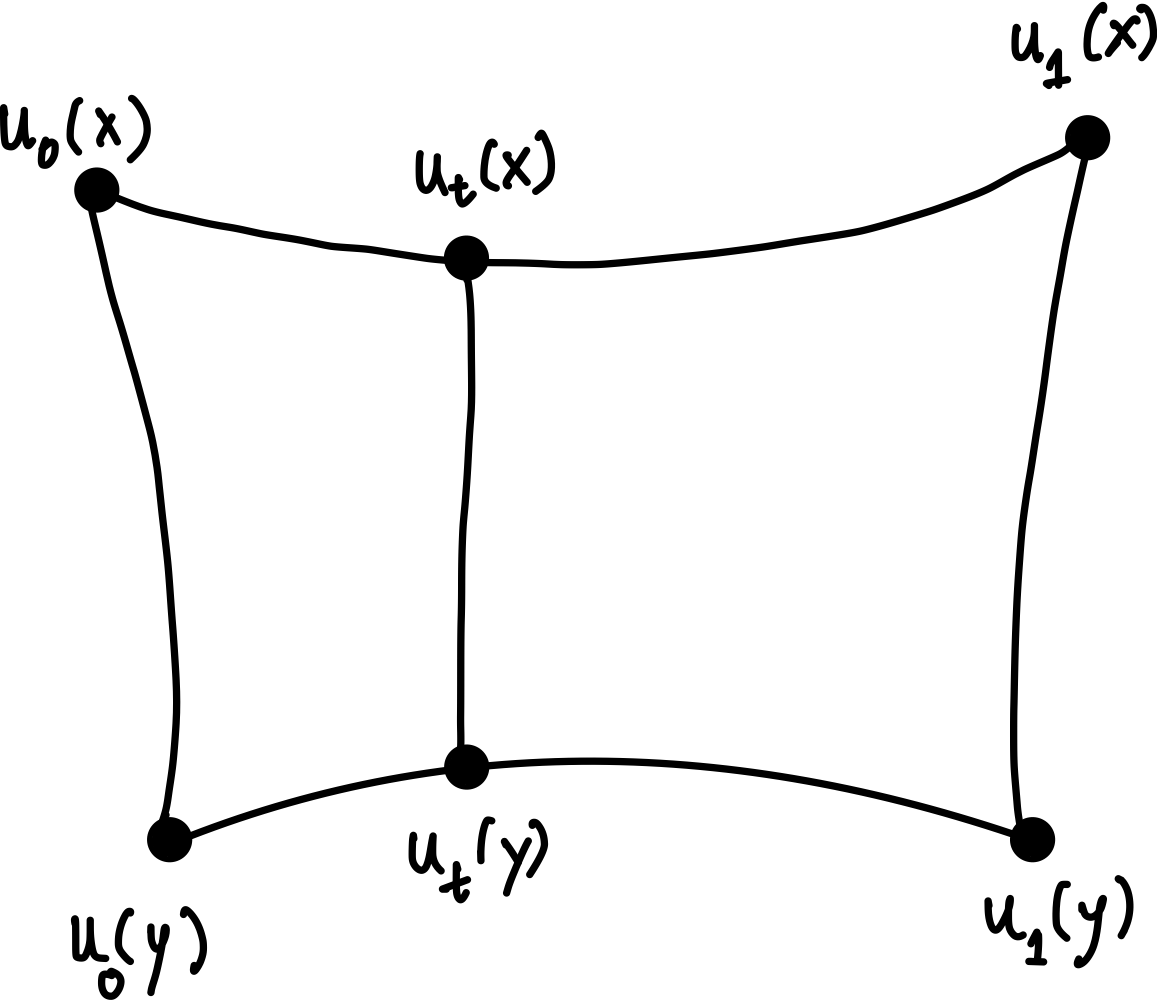}
\caption{$u_0$, $u_1$ and interpolation map $u_t$}%date unknown
\label{quadrilateral}
\end{center}
\end{figure}

This idea leads to the following uniqueness results.

\begin{theorem}[Uniqueness in CAT($\kappa$) with $\kappa<0$ \cite{meseUniqueness}]
If   $X$ has negative curvature (i.e.~every point of $X$  has a neighborhood that is CAT$(\kappa)$ for $\kappa<0$), then the solution to the homotopy or equivariant problem is unique.  
\end{theorem}

\begin{theorem}[Uniqueness in Euclidean buildings \cite{daskal-meseMRL}]
Let $M$ be a finite-volume Riemannian manifold, $X$ an irreducible Euclidean building, and $\rho:\pi_1(M) \rightarrow \mathsf{Isom}(X)$ a homomorphism. If $\rho(\pi_1(M))$  does not preserve any non-empty, closed, convex, proper subset of $X$, then a finite energy $\rho$-equivariant harmonic map $u:\widetilde M \rightarrow X$ is unique.
\end{theorem}

\subsection{Rigidity according to Gromov-Schoen}
\label{sec:GSrigidity}
%\subsection{Motivating examples for regularity theory} \label{sec:motivating}
In \cite{gromov-schoen}, Gromov and Schoen proved the following rigidity result which extends the celebrated Margulis superrigidity theorem \cite{margulis} and Corlette's rank 1 Archimedean superrigidity theorem \cite{corlette}.
  
  \begin{theorem}[Rank 1 $p$-adic superrigidity] \label{thm:rank 1}
  \label{GSrigiditythm}
Let $G / K$ be an irreducible rank 1 symmetric space of noncompact type other than the real and complex hyperbolic space.  Let  $\Gamma$ be a lattice in  $G$ and   
% b; $G / K$ be an irreducible symmetric space of noncompact type other than $SO_0(p,1) \slash SO(p) \times SO(1)$, $SU_0(p,1) \slash S(U(p) \times U(1))$.
  %Let   $\Gamma$ be a lattice in  $G$  and 
  let $\rho: \Gamma \rightarrow \mathsf{Isom}(X)$ be a homomorphism where $X$ is a locally finite Euclidean building.  Then any finite energy $\rho$-equivariant harmonic map $u:G/K \rightarrow X$ is a constant.  
% If  $\Gamma$ is a lattice in a simple rank-1 Lie group excluding the real and complex hyperbolic cases, then every $p$-adic representation $\rho:\Gamma \to \mathrm{GL}(n,\mathbb Q_p)$ for $n \in \mathbb N$ and prime $p$ has a precompact image in $\mathrm{GL}(n,\mathbb Q_p)$.
 \end{theorem}

 Theorem~\ref{GSrigiditythm} together with the equivariance of $u$ implies that $\rho(\pi_1(M))$ fixes a point in $X$.
The significance of the Corlette-Gromov-Schoen result is that it provides the first  superrigidity theorem for rank-1 lattices.
Gromov and Schoen's work was pioneering in bringing harmonic map methods into rigidity theory for singular targets, such as Euclidean buildings, opening the door to many subsequent developments.

The key idea behind the rigidity theorems of Gromov and Schoen is to prove that a harmonic map $u:M \to X$ into an NPC space is regular on a big open set of the domain.  More precisely,
a point $p \in \Omega$ is called a \emph{regular point} of a map $u : \Omega \to X$ if there exists $r > 0$ such that the image of the geodesic ball $B_r(p)$ under $u$ is contained in a part of $X$ that is isometric to a smooth manifold.  Otherwise, $p$ is called a \emph{singular point}. We denote the set of regular points and singular points  by $\mathcal R(u)$ and $\mathcal{S}(u)$, respectively.  
The importance of these definitions is that,  in a neighborhood of a regular point, differential geometric methods can be applied.  
However, in order to employ the Bochner method,  one needs to prove the following regularity properties:  
\begin{itemize}
\item[(i)] (Hausdorff dimension of the singular set) $\dim_{\mathcal H}  \mathcal S(u)\leq m-2$ where $\dim M=m$.
\item[(ii)] (Energy decay at the singular point) There exists $c>0$ and $\epsilon>0$ such that $
\int_{B_r(x_0)} |du|^2 d\mathrm{vol}_g \leq cr^{m+\epsilon}
$ for $x_0 \in \mathcal S(u)$.
\end{itemize}

We now illustrate why these properties may be needed by examining the Eells-Sampson Bochner formula  given in Theorem~\ref{BochnerES}.   Recall that the main idea is to integrate $\Delta |du|^2$  and take advantage of the fact that the integral is zero.  In the presence of singular points this needs further justification:
multiply the Bochner formula by a cut-off function $\varphi$ that vanishes in a neighborhood of the singular set $\mathcal S(u)$.  Assuming that the domain $M$ is compact, 
apply the divergence theorem to obtain
\begin{equation} \label{integrateES}
0 =  \int_M  \mathrm{div} (\varphi  \nabla |du|^2 ) \,d\mathrm{vol}_g = \int_M \varphi \Delta |du|^2\, d\mathrm{vol}_g + \int_M (\Delta \varphi)  |du|^2\, d\mathrm{vol}_g.
\end{equation}

 From regularity property (i),  for any $\delta>0$, there exists a cover of $\mathcal S(u)$ by  finitely many  balls $\{B_{r_i}(x_i)\}_{i=1}^I$ such that $x_i \in \mathcal S(u)$ and $\sum_{i=1}^I r_i^{m-2+\epsilon}<\delta$. Moreover, the cut-off function $\varphi$  satisfies $\sup_{x \in B_{r_i}(x_i)} |\Delta \varphi| \leq c'r_i^{-2}$  for some  constant $c'>0$ with support in $\bigcup_{i=1}^I B_{r_i}(x_i)$. Thus, the Lipschitz continuity of $u$ (cf.~Theorem~\ref{2.4.6}) and regularity property (ii) imply that the second term on the right hand side of (\ref{integrateES}) can be bounded: \[
\left| \int_M (\Delta \varphi) |du|^2 d\mathrm{vol}_g \right|\leq  \sum_{i=1}^I \left| \int_{B_{r_i}(x_i)} (\Delta \varphi) |du|^2 d\mathrm{vol}_g \right| \leq Lc'c \sum_{i=1}^I r_i^{m-2+\epsilon} \leq Lc'c\delta.
 \]
Thus, applying the Eells-Sampson Bochner formula to the first term on the right hand side of (\ref{integrateES}) and assuming non-negative Ricci curvature for $M$ and non-positive sectional curvature for $N$, we obtain
\[
\int_M \varphi |\nabla du|^2 \leq Lcc'\delta.
\]
Since $\delta>0$ is arbitrary, 
we conclude that  $|\nabla du|^2=0$; hence $u$ is totally geodesic.

We next recall examples from \cite[Introduction]{gromov-schoen} for two contrasting scenarios for the size of the singular set. The first gives a counterexample to the general expectation that the singular set of a harmonic map into an NPC space has Hausdorff codimension at least 2, while the second presents a case where the Hausdorff codimension of the singular set  is  2.
  \\

\begin{figure}[h]
\begin{center}
\includegraphics[width=0.6\textwidth]{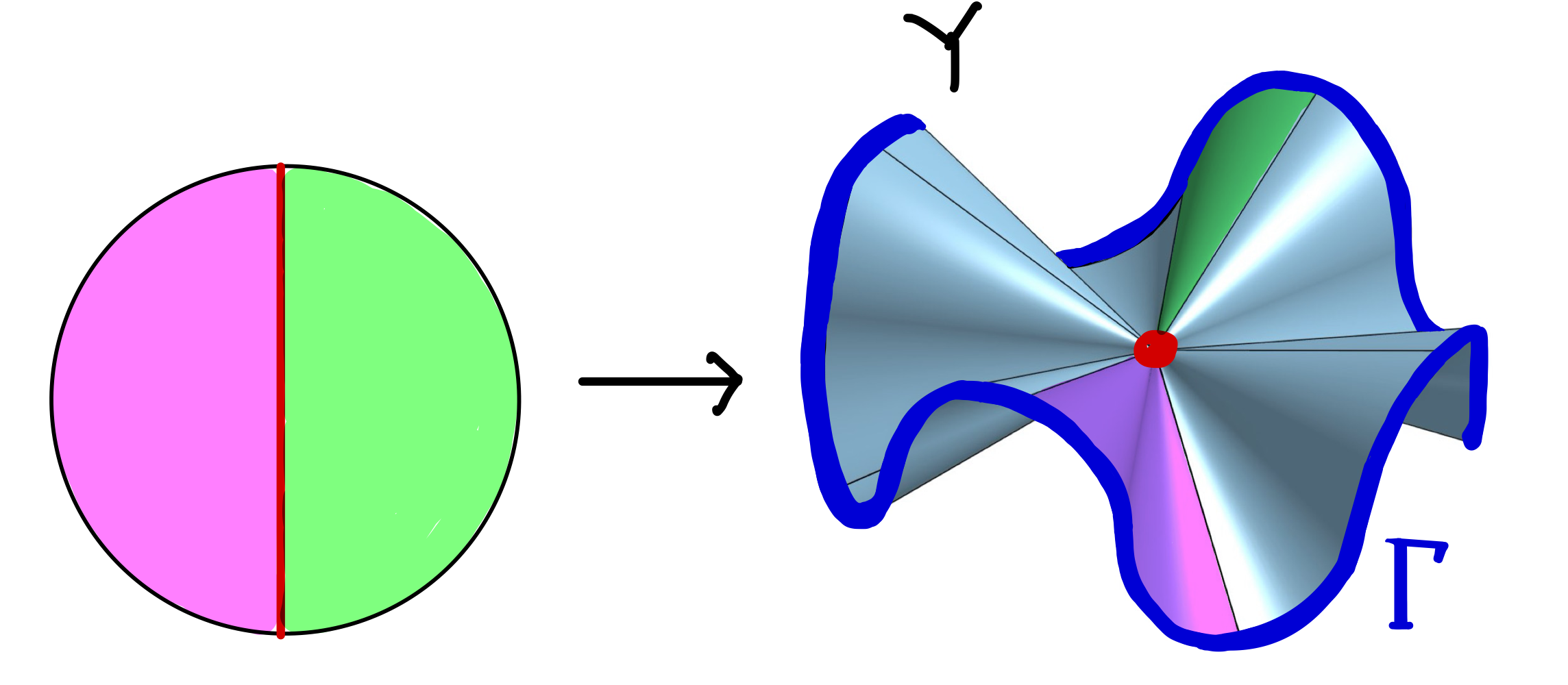}
\caption{Codimension 1 singular set}
\label{codim1}
\end{center}
\end{figure}

{\it Example 1}. This example shows that a harmonic map into an NPC space may have a hypersurface as a singular set:
Let  $Y$ be the unit cone over the circle  of length $\ell > 2\pi$ and  
$M = M_{0} \times [-1,1]$, where $M_{0}$ is a closed manifold. Define a map from the boundary of $M$ 
into  $S^1 \subset  Y$ by sending $M_{0} \times \{-1\}$ to an arc of length less than
$ ( \ell - \pi)/2$, and  $M_{0} \times \{1\}$ to  the symmetrically opposite arc of $S^{1}$ (see Figure~\ref{codim1}).
In $Y$, the convex hull of the two arcs is precisely the union of the cones spanned by them, 
and  the harmonic map $M \to Y$ solving the Dirichlet problem is constrained to lie inside this union.  
The singular set, i.e the set mapped to the  cone vertex, contains a hypersurface separating $M_{0} \times \{-1\}$ from $M_{0} \times \{1\}$.

{\it Example 2.}
This example shows that the singular set of a harmonic map is of codimension 2, given additional structure on the target space. Consider the leaf space of the vertical foliation of the quadratic differential $zdz^2$ on ${\mathbb C}$, endowed with the distance function defined by the transverse measure. This space is isometric to the tripod 
shown in Figure~\ref{tripod}.
  
\begin{figure}[h]
\begin{center}
\includegraphics[width=0.4\textwidth]{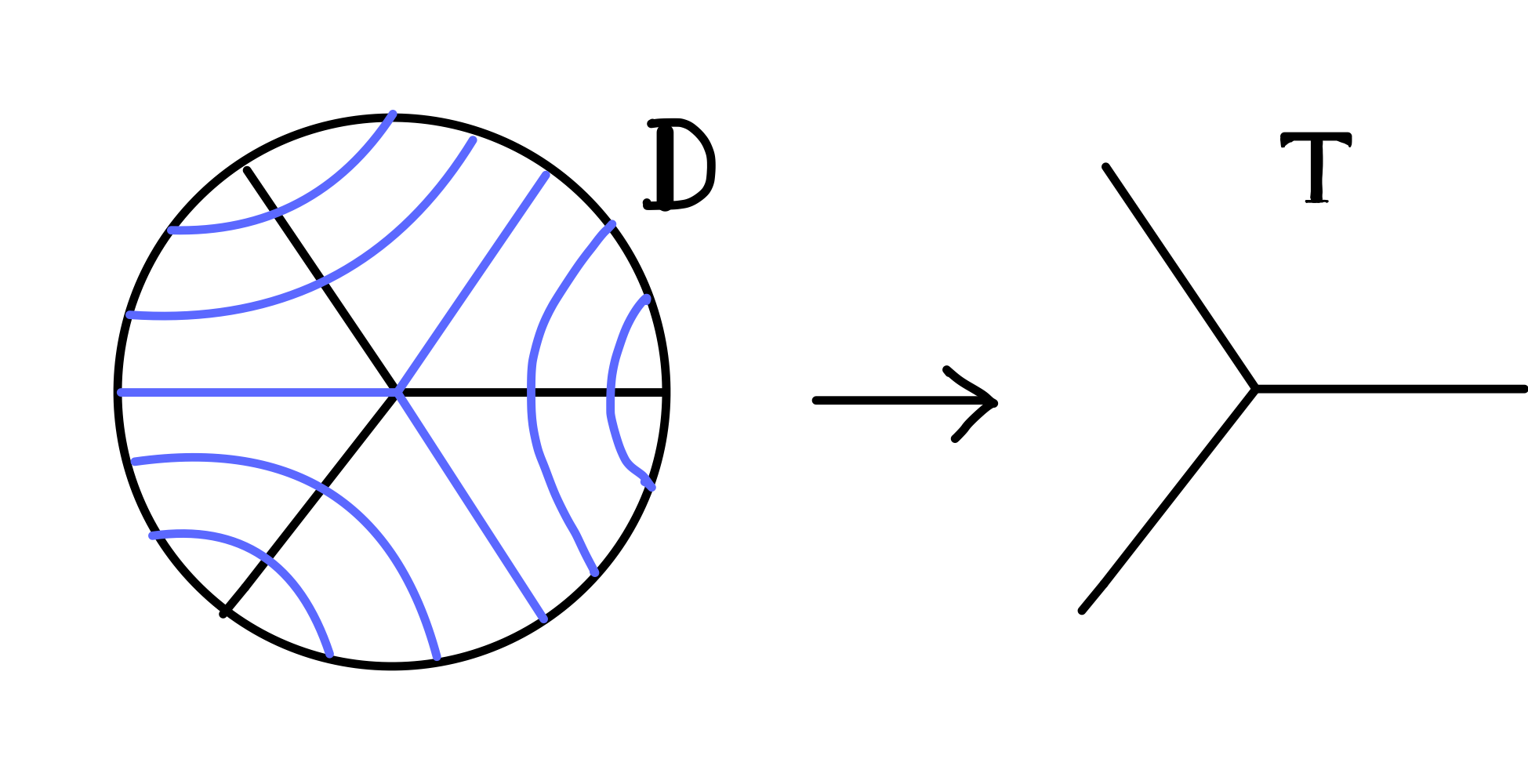}
\caption{Projection of vertical leaves}
\label{leaf}
\end{center}
\end{figure}

  The natural projection map $u: {\mathbb D} \rightarrow T$ from the unit disk centered at 0 is an  energy minimizing  map (see Figure~\ref{leaf}).  It  takes any neighborhood away from the origin into a subset isometric to an open interval. 
 Thus,  $u$ is locally a harmonic function away from the origin.    The image of 
 any neighborhood of the origin is not a manifold. Consequently,  the ``regular set '' ${\mathcal R}(u)={\mathbb D} \setminus \{0\}$ and the ``singular set'' ${\mathcal S}(u)=\{0\}$. 
 
\vspace*{0.12in}
The important distinction between the two spaces in the examples above is the following: for the target $T$ in Example 2, any two points are contained in an isometrically embedded, totally geodesic differentiable manifold (more specifically, in an isometric copy of $\mathbb{R}$).    In Example 1,  for two points at a distance $>\pi$ on the circle of length $\ell > 2\pi$ embedded in cone $Y$, there does not exist any isometrically embedded, totally geodesic differentiable surface containing them.  Given additional structure much like that of Example 2, it is possible to prove the required regularity theorems which allows us to generalize Bochner methods.
We shall see in Section 2.1 that $T$
belongs to a general class of NPC spaces called DM-complexes.

%In the following sections, we will discuss the regularity theory and outline the main ideas of the proofs. This is essential for applying the Bochner method (discussed for smooth targets in Section~\ref{smoothharma}) in the singular setting.
%Recall that the Bochner formula (cf. Theorem~\ref{BochnerES}, or \cite{siu}) is a differential identity involving higher derivatives of the harmonic map, and depends fundamentally on the smooth structure of the Riemannian manifolds involved. In the singular setting 
%In the final step  the Bochner method involves integration by parts, which in the non-smooth setting requires the justification. 
%For this one must first establish that harmonic maps  are  regular on a large open subset of the domain. 
% This result is reliant on the additional structure on the target NPC space.

\section{Rigidity beyond Gromov-Schoen}
\label{beyondGS}
%The idea behind the proofs of the rigidity theorems in Section~\ref{beyondGS} is to establish that harmonic maps are regular on a large open subset of the domain. The local analytic properties of harmonic maps are shown to be tightly controlled by the geometric features of the target NPC space. The resulting regularity theorems follow from a careful analysis of these properties, and the rigidity results presented in this section then follow as a direct consequence.

Mostow established rigidity theorems for symmetric spaces, and Margulis proved superrigidity for symmetric spaces and their analogues, Euclidean buildings. The rigidity theorem of Gromov–Schoen builds on the observation that Euclidean buildings have enough Euclidean structure   to allow the application of analytic techniques. In our work, we expand on this idea: if an NPC space possesses enough structure resembling that of a Riemannian manifold of non-positive sectional curvature, then  analytic techniques can be applied to  this more general setting.

\subsection{Differentiable manifold complexes}

\label{dmcomplex}
A  large class of NPC spaces that is suitable for our purposes  is what we call  \emph{DM-complexes} (shorthand for differentiable manifold complexes). Roughly speaking, a $k$-dimensional DM-complex is an NPC  space that is also a simplicial complex with the property that the unique geodesic between any two points lies in the image of an isometric, totally geodesic embedding of a differentiable manifold (a ``DM''). This includes Gromov and Schoen's \emph{F-connected} (or flat-connected) complexes-such as Euclidean buildings-as well as other examples like hyperbolic buildings and the completion of Teichm\"uller space.

\begin{definition}
\emph{A cell complex $X$ equipped with a metric $G = \{G^S\}$ is called a \emph{Riemannian complex} if each cell $S$ of $X$ is endowed with a smooth Riemannian metric $G^S$ such that  $G^S$ extends smoothly to the boundary of $S$. Furthermore, if $S'$ is a face of $S$, then the restriction of $G^S$ to $S'$ agrees with $G^{S'}$, and $S'$ is totally geodesic in $S$.}
\end{definition}

%\begin{definition}
%\emph{A $k$-dimensional Riemannian complex $(X, G)$ is called a \emph{DM-complex} if, for any two simplices $S_1$ and $S_2$ in $X$ such that $S_1 \cap S_2 \neq \emptyset$, there exists a $k$-dimensional, $C^\infty$-smooth, complete Riemannian manifold $M$ and an isometric, totally geodesic embedding $J : M \rightarrow X$ such that $S_1 \cup S_2 \subset J(M)$. By abuse of notation, we will often denote $J(M)$ by $M$ and refer to it as a DM.}
%\end{definition}

\begin{definition}
\emph{A $k$-dimensional Riemannian complex $(X,G)$ is said to be a DM-complex if every pair of intersecting cells is contained in a smooth, isometrically embedded, totally geodesic $k$-dimensional Riemannian submanifold $M$ in $X$. We will often refer to $M$ as a DM.
}
\end{definition}

\begin{figure}[h]
\begin{center}
\includegraphics[width=0.4\textwidth]{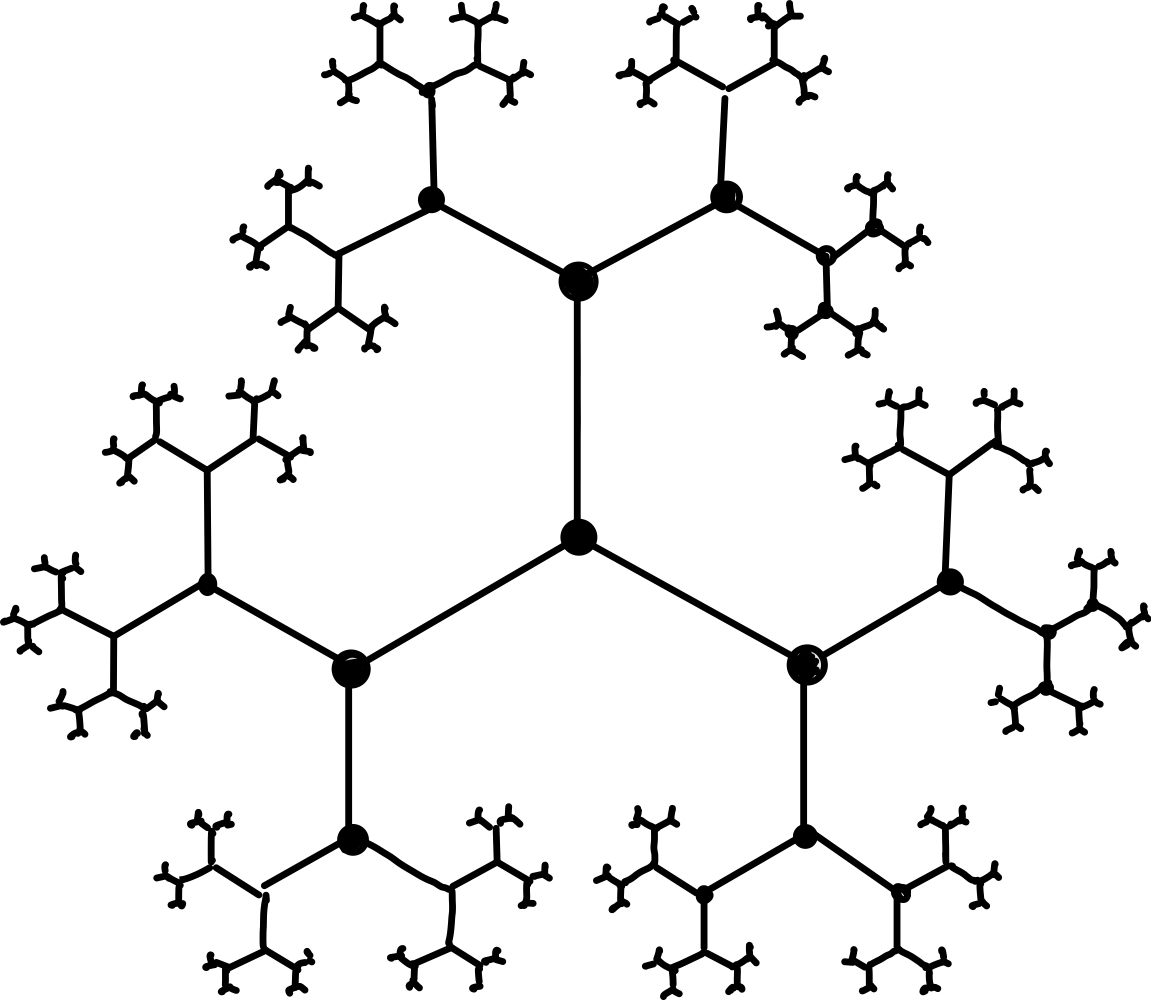}
\caption{1-dimensional Euclidean building}
\label{BTbldg}
\end{center}
\end{figure}

Figure~\ref{BTbldg} is an example of a Euclidean building.  The  distance is represented in a way that is visually distorted near the boundary at infinity, much like in the disk model of hyperbolic space. In particular, every edge has uniform length~$1$ and the local geometry is isometric to the tripod shown in Figure~\ref{tripod}.  Notice that any two points lie in a totally geodesic, isometrically embedded copy of the real line.

The NPC space $Y$ of Example 1 in Section~\ref{sec:GSrigidity} is not a DM-complex.  There does not exist any isometric, totally geodesic embedding of a 2-dimensional smooth Riemannian manifold that contains the cone point.

\subsection{Hyperbolic buildings}

\label{hyperbolicbuildings}
Hyperbolic buildings are analogues of Euclidean buildings in which the apartments are hyperbolic spaces. More precisely, a hyperbolic building is a DM-complex formed from simplices in hyperbolic space $\mathbb{H}^k$ with the property that for any two adjacent cells there exists an isometric and totally geodesic embedding $J : \mathbb{H}^k \to X$ whose image $J(\mathbb{H}^k)$ contains both cells. In this section we assume that all buildings are locally compact.

To address the rigidity of harmonic maps to hyperbolic buildings, we introduce the following:
    
     \begin{definition} Let $\Gamma$ be a discrete group, $X$ an NPC space and $\mathsf{Isom}(X)$ the group of isometries of $X$.   A homomorphism $\rho: \Gamma \rightarrow \mathsf{Isom}(X)$  is called \emph{reduced} if 
(i) $\rho$ does not fix a point at infinity of $X$; and
(ii) there is no unbounded closed convex  $Z \subset X$,  $Z \neq X$ preserved by
 $\rho(\Gamma)$. 
 \end{definition}

 The above definition is the generalization in the NPC setting of the condition that $\rho$ has Zariski dense image in Margulis' theorem.
%\begin{theorem}[Geometric rigidity \cite{daskal-meseSR}] \label{hyperbolicrigidity}
% Let $G / K$ be an irreducible symmetric space of noncompact type other than the real and complex hyperbolic spaces.
%  Let 
% $\Gamma$ be a lattice in  $G$  and let $\rho: \Gamma \rightarrow \mathsf{Isom}(X)$ be a  group homomorphism where $X$ is an NPC DM-complex. 
% %{\color{red} I DON'T THINK WE NEED THIS STATEMENT BECAUSE THE ASSERTION IS ABOUT HARMONIC MAPS.  THIS STATEMENT IS FOR EXISTENCE.  If the rank of $G/K$ is $\geq 2$,  we assume additionally that $\Gamma$ is cocompact.} 
% If the rank of $G/K$ is 1, we assume additionally that the curvature operator of any DM in $X$ is non-positive.  Then any finite energy $\rho$-equivariant harmonic map $u:G/K\rightarrow X$ is non-branching and totally geodesic.   In other words, the image of $u$ is contained in a single DM  of $X$ and  $u$ is totally geodesic as a map into that DM.  
% \end{theorem}

\begin{theorem}[Geometric rigidity in hypebolic buildings \cite{daskal-meseSR}] \label{hyperbolicrigidity}
 Let $G / K$ be an irreducible symmetric space of noncompact type other than the real and complex hyperbolic spaces, $\Gamma$ be a lattice in $G$, and $X$ be an NPC DM-complex.
Assume that either
\begin{itemize}
\item[(a)] $\mathrm{rank}(G/K) \geq 2$ and $\Gamma \backslash G/K$ is compact, or
\item[(b)] $\mathrm{rank}(G/K) = 1$.
\end{itemize}
Then, for any reduced homomorphism $\rho: \Gamma \rightarrow \mathsf{Isom}(X)$,  a finite energy $\rho$-equivariant harmonic map $u:G/K\rightarrow X$ is non-branching and totally geodesic.   In other words, the image of $u$ is contained in a single copy of hyperbolic space isometrically and totally geodesically embedded in  $X$ and  $u$ is totally geodesic as a map into that hyperbolic space.  
 \end{theorem}

  A  consequence of combining the existence theorem of Theorem~\ref{existence} and the geometric rigidity result of Theorem~\ref{hyperbolicrigidity}  is the  superrigidity  result for hyperbolic buildings.

 \begin{theorem}[Superrigidity \cite{daskal-meseSR}]  \label{hyperbolicmargulis}
 Let $G \slash K$ be an irreducible symmetric space of noncompact type, other than the real and complex hyperbolic spaces.
 %$SO_0(p,1) \slash SO(p) \times SO(1)$, $SU_0(p,1) \slash S(U(p) \times U(1))$.
  Let 
 $\Gamma$ be a lattice in $G$  and let $\rho: \Gamma \rightarrow \mathsf{Isom}(X)$ be a reduced homomorphism where $X$ is a  hyperbolic building. If the rank of $G/K$ is $\geq 2$,  assume additionally that $\Gamma$ is cocompact. 
 Then $\rho(\Gamma)$ fixes a point in $X$.
  \end{theorem}

Next, recall that a harmonic map from a K\"ahler manifold to a Riemannian complex is called {\it pluriharmonic} if it is pluriharmonic in the usual sense on the regular set. In other words, if its restriction to any embedded holomorphic disc is  harmonic. The following holds:

\begin{theorem}[Pluriharmonicity \cite{daskal-meseSR}] \label{hyperbolicph} Let $\widetilde M$ be the universal cover of a complete finite volume K\"ahler manifold $(M, \omega)$. Let 
$\Gamma =\pi_1(M)$, $X$  a hyperbolic building and $\rho:\Gamma \rightarrow \mathsf{Isom}(X)$ a homomorphism. 
Then any finite energy $\rho$-equivariant harmonic map $u: \widetilde M \rightarrow X$ is pluriharmonic.
\end{theorem}

Theorems~\ref{hyperbolicmargulis} and~\ref{hyperbolicph} mirror the corresponding theorems for Euclidean buildings (cf. \cite{gromov-schoen}).

\subsection{Teichm\"uller space}
\label{teichmuller}
%\subsection{The Weil-Petersson metric on Teichm\"uller space} 
 Let $\mathcal T$ denote the Teichm\"{u}ller space  parametrizing surfaces  $S$ of   genus $g$  with $p$-marked points  and $3g-3+p>0$.  For $\sigma \in \mathcal T$ represented (by a slight abuse of notation) as a complete hyperbolic metric $ds^2=\sigma |dz|^2$,  the cotangent space $T^*_\sigma(\mathcal T)$ can be identified with the space of meromorphic quadratic differentials $\Phi=\phi\, dz^2$ on $S$ with respect to $\sigma$, having at most simple poles at the marked points. The Weil-Petersson metric $G_{WP}$  is defined by
\[
||\Phi||^2_{G_{WP}}= \int_S |\phi(z)|^2\sigma^{-1}dxdy
\]
 and the space $(\mathcal T, G_{WP})$ is an incomplete  K\"ahler  manifold (cf.~\cite{wolpertPJ}, \cite{chu}) of negative sectional curvature (cf.~\cite{tromba} and \cite{wolpertJDG}). Moreover, the curvature tensor of $G_{WP}$ is strongly negative  in the sense of Siu \cite{schumacher}.  This suggests utilizing Siu's Bochner formula in order to study the holomorphic rigidity of Teichm\"uller space.  
 
 Since the Weil-Petersson metric is incomplete, we first need to take the completion (equivalence classes of Cauchy sequences in $\mathcal T$) denoted by $\overline{\mathcal T}$. This is nothing but the {\it augmented Teichm\"uller space}
 (cf.~ \cite{abikoff}). The boundary $\partial {\mathcal T}=\overline{\mathcal T} \setminus \mathcal T$  is stratified by (totally geodesic) disjoint smooth strata $\mathcal T_c$, each parametrized by collections  $c$ of disjoint simple closed essential and non-peripheral curves.  More precisely,
 given $c$,  $\mathcal T_{c}$ denotes the Teichm\"uller space of the
   nodal surface obtained by collapsing each of the curves in $c$ to a point.   Let $d_{\overline{\mathcal T}}$ denote the distance function induced by $G_{WP}$. The space    
 $(\overline{\mathcal T},d_{\overline{\mathcal T}})$,  is an NPC space   (cf.~ \cite{yamada}, \cite{daskal-wentworth}).

A large part of the harmonic map theory for NPC targets discussed earlier applies to  $({\mathcal T},d_{\overline{\mathcal T}})$.     
However, there are some major difficulties: The most significant one is that  $\overline{\mathcal T}$   is not locally compact. This is in sharp contrast to the Euclidean buildings studied by Gromov and Schoen \cite{gromov-schoen} and the hyperbolic buildings discussed in Section~\ref{hyperbolicbuildings}. This poses a serious difficulty in the study of harmonic maps, which we will explain further in the next sections. For now, we state the main rigidity result for   Teichm\"uller space:

\begin{theorem}[Holomorphic rigidity \cite{daskal-meseHR}]\label{rigiditytheorem0}
Let $M$ be a complete, finite volume K\"ahler manifold  with universal cover $\widetilde M$ and $\pi_1(M)$ finitely generated. Let $ \Gamma$ be the mapping class group of an oriented  surface  $S$ of genus $g$  and  $p$ marked points such that $k=3g-3+p>0$, $\overline {\mathcal T}$  the Weil-Petersson completion of the Teichm\"{u}ller space ${\mathcal T}$ of $S$ and $\rho: \pi_1(M) \rightarrow \Gamma$ a homomorphism.   If there exists a  finite energy
 $\rho$-equivariant  harmonic map $ u: \widetilde M \rightarrow \overline {\mathcal T}$,  
then there  exists a unique stratum ${\mathcal T}_c$ of $\overline{\mathcal T}$ such that $ u$ defines a pluriharmonic map into ${\mathcal T}_c$.
% Furthermore,
%\[
% \sum_{i,j,k,l}R_{ijkl}d''u_i \wedge d'u_j \wedge d'u_k \wedge d''u_l \equiv 0
%\]
%where $R_{ijkl}$ denotes the Weil-Petersson curvature tensor. In particular,
If additionally the (real) rank of $u$ is $\geq  3$ at some point, then  $ u$ is holomorphic or conjugate holomorphic. 
\end{theorem}

The assumption that a  finite energy
 $\rho$-equivariant  harmonic map to the Weil-Petersson completion $\overline{\mathcal T}$ of Teichm\"{u}ller space exists holds in many important cases, for example when $M$ is closed.    In Section~\ref{sec: infinite energy} we will also discuss rigidity of maps of infinite energy.
 
   Recall from \cite[p.142]{papa} or \cite[Definition 2.1]{daskal-wentworth}  that   two
  pseudo-Anosov elements of the mapping class group are called {\it independent} if their fixed point sets in the space of projective measured foliations  do not coincide.
 A subgroup of the mapping class group $\Gamma$ is called {\it {sufficiently large}} if it contains two independent pseudo-Anosov elements. A homomorphism $\rho$ into the mapping class group is called sufficiently large if its image is sufficiently large. 
% Also recall the notion of a proper action on an NPC space. Let $\rho$ be a homomorphism into the isometry group of an NPC space $(X, d)$. We call $\rho$  {\it proper}
% if the sublevel sets of the function 
%\[
%\delta(P)=\max \{d(\rho(\lambda)P,P): \lambda \in \Lambda\}
%\]  
%are bounded in $\widetilde X$; i.e.~given arbitrary $c>0$, there exists $P_0 \in X$ and $R_0>0$ such that
%\[
%\{P \in \widetilde X:  \delta(P)\leq c\} \subset B_{R_0}(P_0).
%\]
% According to \cite{daskal-wentworth}, if \( \rho \) is sufficiently large, then  it is proper; that is, \( \rho \) satisfies the assumption of Theorem~\ref{rigiditytheorem0}.  Combining with  \cite[Corollary 1.3]{daskal-wentworth}, 
 According to \cite[Corollary 1.3]{daskal-wentworth}, sufficiently large 
$\rho$
guarantees the existence of a finite-energy harmonic map. Combining this with Theorem~\ref{rigiditytheorem0},

   \begin{theorem}\label{rigiditytheorem1}Let $M$ be a complete, finite volume K\"ahler manifold  with universal cover $\widetilde M$ and $\pi_1(M)$  finitely generated. Let $ \Gamma$ be the mapping class group of an oriented  surface  $S$ of genus $g$  and  $p$ marked points such that $k=3g-3+p>0$ and $\rho: \pi_1(M) \rightarrow \Gamma$ a homomorphism that is sufficiently large.   If there exists a  finite energy
 $\rho$-equivariant   map $ \widetilde M \rightarrow  {\mathcal T}$,  
then there  exists a $\rho$-equivariant  pluriharmonic map $ u: \widetilde M \rightarrow  {\mathcal T}$.
 Furthermore,
\[
 \sum_{i,j,k,l}R_{ijkl}d''u_i \wedge d'u_j \wedge d'u_k \wedge d''u_l \equiv 0
\]
where $R_{ijkl}$ denotes the Weil-Petersson curvature tensor. In particular,
if additionally the (real) rank of $u$ is $\geq  3$ at some point, then  $ u$ is holomorphic or conjugate holomorphic. 
\end{theorem}

   The rank condition  holds in many important applications and is usually verified by showing that certain nontrivial homology classes in $M$ of degree $\geq 3$ are mapped  nontrivially under $u$ (see for example \cite{siu}).

\section{Regularity results for DM-complexes}

\subsection{Asymptotic product structure}
 Gromov and Schoen's regularity theorem \cite[Theorem~6.3]{gromov-schoen}  uses  an inductive argument on the dimension of  the Euclidean building.   The core step, \cite[Theorem~5.1]{gromov-schoen}, shows that near a point of order 1 (see definition below), the image of a harmonic map lies in the product of a Euclidean space and a building of lower dimension. Thus, the  harmonic map locally decomposes into two harmonic maps.  More precisely,
 a harmonic map map $u : \Omega \rightarrow X^k$ into a locally finite $k$-dimensional Euclidean building can be written  near a singular point $x_0 \in \mathcal{S}(u)$ of order 1, as 
\begin{equation} \label{Vv}
u = (V, v)
\end{equation}
where $V$ is the nonsingular component mapping into a Euclidean factor $\mathbb{R}^j$, and $v$  mapping into a lower-dimensional complex $X_2^{k-j}$. Both maps $V$ and $v$ are harmonic.  
%The singular set is decomposed as 
%$$
%\mathcal S(u)=\bigcup_j \mathcal S_j(u)
%$$
%where $x \in \mathcal S_j(u)$ means that $u$ has the decomposition (\ref{Vv}) and $x$ is a singular point of the map $v$.

A straightforward application of Federer's dimensional reduction implies that the set of points of order greater than 1 has Hausdorff codimension  2.  Therefore, the analysis of the singular set of $u$ may be inductively reduced to the study of $v$ near  order 1 singular points.  Recall that $v$ maps into the building $X_2^{k-j}$, where $j>0$. However, this inductive step breaks down if the spaces lack such a product structure, making it harder to apply  to a broader class of NPC spaces such as hyperbolic buildings  or the completion of Teichm\"uller space.
%To emphasize this, recall that a harmonic map $u = (u^1, \dots, u^m) : \Omega \rightarrow \mathbb{R}^m$ into Euclidean space can be viewed as $m$ independent harmonic functions. Assuming continuity, a harmonic map between Riemannian manifolds can also be locally expressed in coordinates as a collection of component functions $u = (u^1, \dots, u^m)$. However, when the target metric is non-Euclidean, the component functions are no longer independent. Indeed, the harmonic map equations
%(\ref{hme}) intertwine the components through the Christoffel symbols $\Gamma^i_{jk}$ of the target metric. On the other hand, Riemannian manifolds are locally asymptotic to Euclidean space: in normal coordinates centered at a point, the manifold is Euclidean up to second order. In this sense, a smooth Riemannian manifold is asymptotically a product of $m$ copies of $\mathbb{R}$.

\subsection{Statement of the theorems}

Paper \cite{daskal-meseDM} initiated the  study of harmonic maps into spaces that are asymptotically a product of NPC spaces. In many respects, this transition resembles the passage from harmonic functions to harmonic maps into Riemannian manifolds.
Locally, the DM-complexes discussed in Section~\ref{dmcomplex} have an asymptotic product structure. Thus, even though the maps $V$ and $v$ in the decomposition~\ref{Vv} are not harmonic, they are nevertheless asymptotically harmonic. The precise assumptions on DM-complexes  required to establish the two regularity  theorems below  are detailed in \cite[Chapter~5]{daskal-meseDM}, and we will not discuss them here. 
%These assumptions are satisfied by important examples such as hyperbolic buildings  and Teichm\"uller space equipped with the Weil-Petersson metric.

%The following two theorems constitute the main regularity theorems for DM-complexes. 

\begin{theorem}[Regularity theorem I \cite{daskal-meseDM}] \label{regtheorem} 
Let $u : \Omega \rightarrow X$ be a harmonic map from an $m$-dimensional Riemannian domain into a $k$-dimensional NPC DM-complex. Then the singular set $\mathcal{S}(u)$ has Hausdorff codimension at least 2 in $\Omega$; that is,
\[
\dim_{\mathcal{H}}(\mathcal{S}(u)) \leq m - 2.
\]
\end{theorem}

The next result quantifies the decay of the energy of a harmonic map near its singular points.

\begin{theorem}[Regularity theorem II \cite{daskal-meseDM}] \label{goto0DM}
Let $u : \Omega \rightarrow X$ be as in Theorem~\ref{regtheorem}. For any compact subdomain $\Omega_1 \subset \Omega$, there exists a sequence of smooth functions $\{\psi_i\}$ such that $\psi_i \equiv 0$ in a neighborhood of $\mathcal{S}(u) \cap \overline{\Omega_1}$, $0 \leq \psi_i \leq 1$, and $\psi_i(x) \rightarrow 1$ for all $x \in \Omega_1 \setminus \mathcal{S}(u)$, with
\[
\lim_{i \rightarrow \infty} \int_{\Omega} |\nabla \nabla u| |\nabla \psi_i| \, d\mu = 0.
\]
\end{theorem}

Since Euclidean buildings are special cases of DM-complexes, the regularity results of \cite{gromov-schoen} can be viewed as special cases of the above theorems. Hyperbolic buildings are also special cases of DM-complexes, hence Theorems~\ref{regtheorem} and~\ref{goto0DM} apply. Combined with the appropriate Bochner formula, they imply Theorems~\ref{hyperbolicrigidity} and~\ref{hyperbolicmargulis}. 

\subsection{Monotonicity formulae}

We briefly discuss the main tools used in the proof of Theorems~\ref{regtheorem} and~\ref{goto0DM}.  The first  is the monotonicity formula for harmonic maps.  
Studied by Gromov-Schoen in the context of NPC spaces, this formula originated in \cite{almgren} and is often referred to as the {\it frequency function}.   In \cite{daskal-meseDM}, we modified the formula to apply to maps that are only asymptotically harmonic by adding correction terms. This implies
%  We thus obtain the monotonicity inequality
%\[
%\frac{1}{r} + \frac{E'(r)}{E(r)} - \frac{I'(r)}{I(r)} \geq -C,
%\]
%where $x_0 \in \mathcal S(u)$,  
%\begin{eqnarray*}
%E(r) & = & \int_{B_r(x_0)} |\nabla v|^2 d\mathrm{vol}_g
%\\
%I(r) & = & \int_{\partial B_r(x_0)} d(v,v(0)) d\Sigma_g
%\end{eqnarray*}
%and  $C$ is a constant depending on the deviation of the domain metric $g$ from te Euclidean metric.  As a consequence, there exists an order function for the singular component map $v$ of $u$ and use it to analyze its behavior.  
\begin{theorem}[Order of the Singular Component \cite{daskal-meseDM}] \label{limitexists*}
Let $u : \Omega \rightarrow X$ be a harmonic map from an $n$-dimensional Riemannian domain into a $k$-dimensional NPC DM-complex. If $x_0$ is a singular point of order one and $u = (V, v)$  near  $x_0$, then
\[
\mathrm{Ord}^{v}(x_0) := \lim_{\sigma \rightarrow 0} \frac{\sigma E(\sigma)}{I(\sigma)}
\]
exists, where $E(\sigma)= \int_{B_\sigma(x_0)} |\nabla v|^2 d\mathrm{vol}_g$ and 
$I(\sigma) = \int_{\partial B_\sigma(x_0)} d(v,v(0)) d\Sigma_g$.
\end{theorem}   
We remind the reader that, for harmonic functions, the order is the degree of the dominant homogeneous harmonic polynomial approximating $u-u(x_0)$ near $x_0$.

%\begin{theorem}[Gap Theorem \cite{daskal-meseDM}] \label{gaptheorem*}
%Under the same assumptions as Theorem~\ref{limitexists*}, there exists $\epsilon_0 > 0$ such that $\mathrm{Ord}^v(x) \geq 1 + \epsilon_0$ for all $x \in \mathcal{S}_j(u)$ near $x_0$.
%\end{theorem}

At  order 1 points,  a harmonic map $u$ can be closely approximated by a homogeneous degree 1 map $l$.  This means that $t \mapsto l(t\xi)$ is a constant speed geodesic from $l(0)$ to $l(\xi)$ for any $\xi \in \partial B_1(0)$, and that there exists a function $\mu(\sigma)$ with $\lim_{\sigma \to 0}\mu(\sigma)= 0$  and
\begin{equation} \label{linap}
\sup_{x \in B_\sigma(0)} d(u(x), l(x)) < \mu(\sigma).
\end{equation}

\subsection{Effectively contained and essentially regular} \label{sec:gromov-schoen}

We next introduce two fundamental concepts  from the work of Gromov and Schoen \cite{gromov-schoen}
which are key in proving the regularity Theorems \ref{regtheorem}  and \ref{goto0DM}. 
Let $X$ be an NPC space
and  $X_0$ a totally geodesic subspace of $X$.  
The first  is the notion  of   a homogeneous degree 1 map  $l: \R^n \rightarrow X_0 \subset X$ being \emph{effectively contained} in $X_0$.  One may interpret this as asserting that a sufficiently small neighborhood of the image of $l$ lies in $X_0$, except for a set of small measure.  
The second notion is the \emph{essential regularity} of $X_0$. This means that any harmonic map into $X_0$ admits an approximation by a homogeneous map  that is better than first order.     
To illustrate these notions, we give the following example, which forms the first step in Gromov-Schoen's inductive proof of regularity.\\

\begin{figure}[h] 
  \includegraphics[width=2.5in]{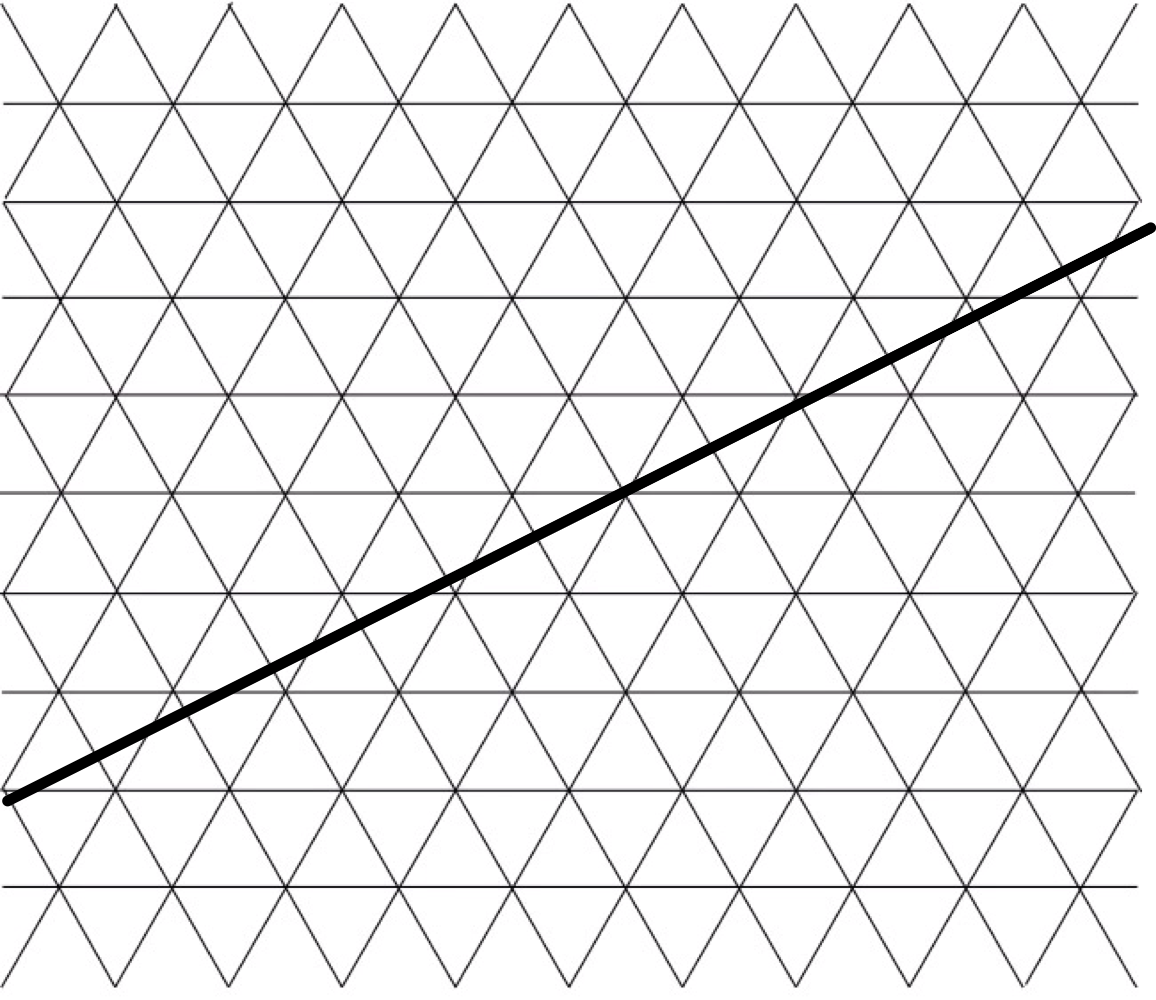}
 \caption{In this diagram, the set of points on the thick line close to the complement of an apartment is small.  }
 \label{diagram:effective}
 \end{figure}

\emph{Example.} Let $X$ be the $k$-pod from Section~\ref{npctarget} and with vertex $V$. (For $k=3$, this is the tripod of Figure~\ref{tripod}.)  Let $X_0 \subset X$ be a union of two copies of $[0,\infty)$.  Thus, $X_0$ is a totally geodesic subspace isometric to $\R$. Let
\begin{equation} \label{alphaxb}
l: \R^n \rightarrow \R \simeq X_0 \subset X, \quad l(x) = \vec{A} \cdot x
\end{equation}
be a linear function for some $\vec{A} \in \R^n$. Here, $0 \in \R \simeq X_0$ is identified with the vertex of $X$. Rotating coordinates if necessary, we may assume $\vec{A} = (a, 0, \ldots, 0)$ so that $l(x) = a x^1$. Let ${\bf B}_{\epsilon}(X_0)$ denote the $\epsilon$-neighborhood of the set $X_0$. Then,
\[
{\bf B}_{\delta \sigma}(l(x)) \cap (X \setminus X_0) \neq \emptyset 
\quad \Leftrightarrow \quad |l(x)| < \delta \sigma 
\quad \Leftrightarrow \quad |x^1| < \frac{\delta}{|a|} \sigma.
\]
Consequently, for any $\epsilon > 0$, there exists $\delta > 0$ (for instance, $\delta = \frac{\epsilon |a|}{2 v_{n-1}}$, where $v_{n-1}$ is the Euclidean volume of the unit $(n-1)$-ball) such that
\begin{equation} \label{exeffct}
\mathrm{Vol} \left\{ x \in B_\sigma(0) : {\bf B}_{\delta \sigma}(l(x)) \cap (X \setminus X_0) \neq \emptyset \right\} < \epsilon \sigma^n.
\end{equation}
For any $x$ sufficiently close to $0$ and $P \in X_0$ sufficiently close to the vertex $V$, define a homogeneous degree 1 map $l_{x,P}$ by requiring that for each $\xi \in \partial B_1(x)$, the curve $t \mapsto l_{x,P}((1-t)x+t\xi)$ is a geodesic from $l_{x,P}(x)=P$ to $l_{x,P}(\xi)=l(\xi)$. Then (\ref{exeffct}) continues to hold with $l$ replaced by $l_{x,P}$.    According to \cite[Definition p.~211]{gromov-schoen},  this is  the formal definition of the linear map $l$ being   \emph{effectively contained} in $X_0$.  (See Figure~\ref{diagram:effective}.)

Moreover,  the totally geodesic subspace $X_0$ is \emph{essentially regular}. For any harmonic function $f: (B_1(0), g) \to X_0 \simeq \R$, the Taylor approximation yields
\[
d(f(x), l(x)) \leq C |x|^2,
\]
where $l(x) = \nabla f(0) \cdot x + f(0)$ and $C$ depends only on the geometry of the domain and the total energy of $f$. Thus, there exist constants $\alpha > 0$ (in this case, $\alpha = 1$) and $C > 0$ such that
\begin{equation} \label{exessreg}
\sup_{x \in B_\sigma(0)} d(f(x), l(x)) \leq C \sigma^{1+\alpha} \sup_{x \in B_1(0)} d(f(x), L(x)), \quad \forall \sigma \in (0, \tfrac{1}{2}],
\end{equation}
for any affine function $L(x) = \vec{A} \cdot x + b$ (cf. \cite[Definition p.~210]{gromov-schoen}).
Key to the definition of essential regularity is that the constants $\alpha$ and $C$ are independent of the subspace $X_0$, depending only on the geometry of the domain and the total energy $E^f$ of $f$.

The notions of \emph{effectively contained} and \emph{essentially regular} extend to higher-dimensional locally finite simplicial complexes, such as buildings.

\subsection{The iteration}
We now illustrate the iterative technique used by Gromov Schoen to prove the regularity theorems.   For simplicity, we focus on a harmonic map $u: B_1(0) \to X$ where $X$ is a $k$-pod with vertex $V$ as in the previous section.  Let $X_0 \approx \mathbb R$ be the union of two copies of $[0,\infty)$.  Assume that $l: B_1(0) \to X_0$ is a homogeneous degree 1 map effectively contained in $X_0$, and that $u$ and $l$ are $D$-close in the sense that
\begin{equation} \label{setup}
\sup_{x \in B_1(0)} d(u(x), l(x)) < D.
\end{equation}
The goal is to produce a \emph{linear scale approximation}:
\begin{equation} \label{ssa}
\sup_{x \in B_\sigma(0)} d(u(x), l(x)) < c \sigma, \quad c > 0.
\end{equation}
In other words, $\mu(\sigma)= c\sigma$ in (\ref{linap}).

\begin{figure}[h]
\begin{center}
\includegraphics[width=0.4\textwidth]{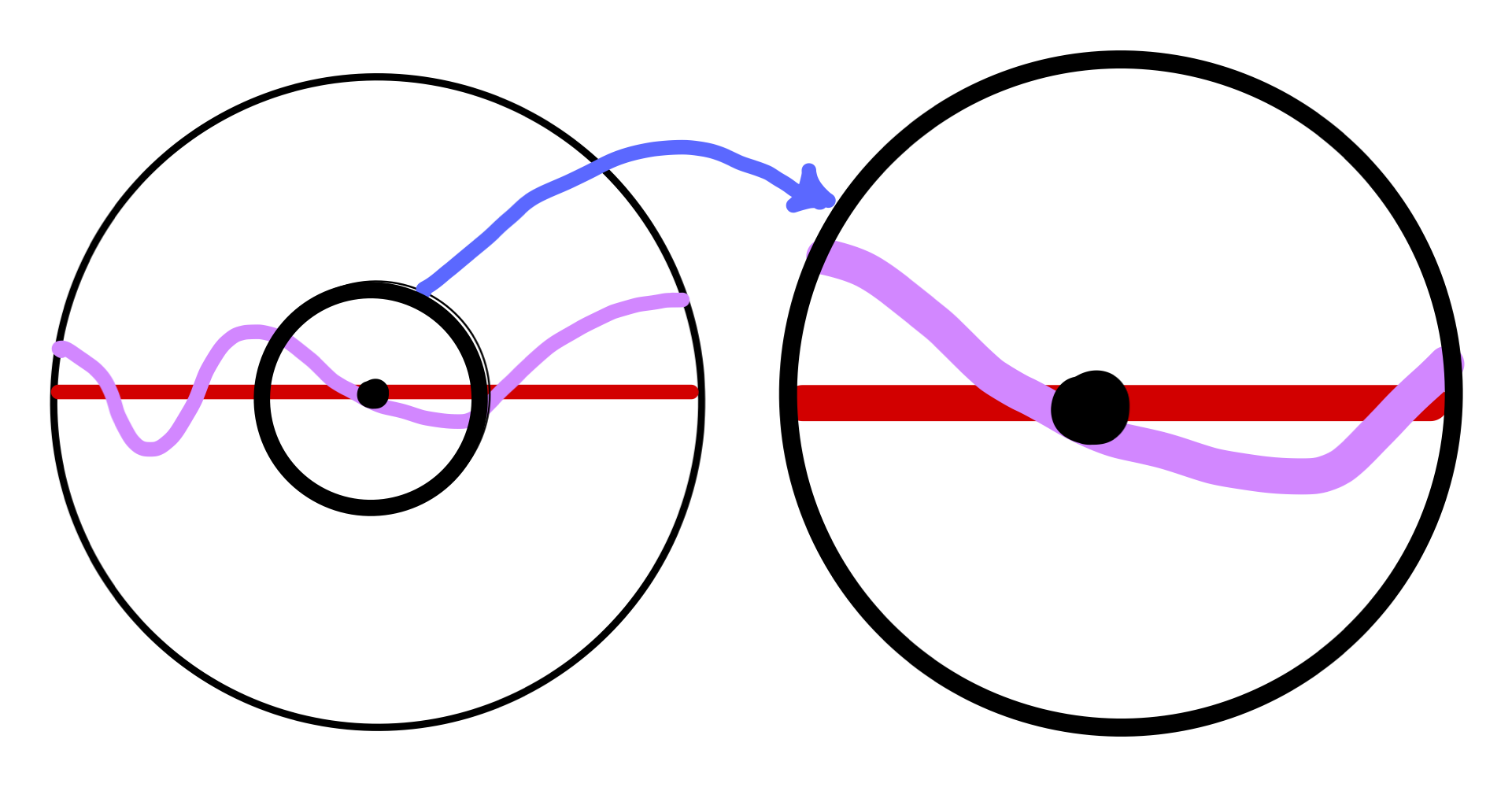}
\caption{Inductive iteration scheme}
\label{blowup}
\end{center}
\end{figure}

Estimate~\eqref{ssa} is obtained inductively through an iteration scheme: at each stage, one improves the approximation by a fixed factor. (See Figure~\ref{blowup}.)
Specifically, there exists $\theta \in (0, \tfrac{1}{2}]$ such that if an affine map
\[
{}_i l: B_{\theta^i}(0) \to X_0
\]
approximates $u$ at the $i$th stage, then there exists a new affine map
\[
{}_{i+1} l: B_{\theta^{i+1}}(0) \to X_0
\]
that approximates $u$ more closely in the smaller ball $B_{\theta^{i+1}}(0)$ at the $(i+1)^{th}$ stage.

To construct $_{i+1}l$, consider the harmonic map $v: B_{\theta^i}(0) \to X_0$ with boundary data $\pi \circ u$, where
\[
\pi: X \to X_0
\]
is the nearest point projection. Since $X_0$ is essentially regular, $v$ admits a good linear approximation $_{i+1}l$. The key technical step relies on the fact that if $l$ is effectively contained in $X_0$ and approximates $u$, the set of points $x \in \partial B_{\theta^{i+1}}(0)$ such that $\pi \circ u(x)\neq u(x)$ is of small measure. Combining this with mean value inequality for the subharmonic function $d(u,v)$ on $B_{\theta^{i+1}}(0)$,  one concludes that $u$ and $v$ are close. Consequently, the affine map $l_{i+1}$ provides a better approximation to $u$ than the previous approximation $l_i$;  indeed,
\[
\sup_{x \in B_{\theta^{i}}(0)} d(u(x),{_{i}l}(x)) \leq \theta^{i}d_0 \ \ \Rightarrow \ \ \sup_{x \in B_{\theta^{i+1}}(0)} d(u(x),{_{i+1}l}(x)) \leq \theta^{i+1}\frac{d_0}{2}.
\]
Combining the above with the triangle inequality and the homogeneity of $l$ and $_i l$, 
\[
\sup_{x \in B_{\theta^{i}}(0)} d(l(x),{_{i}l}(x)) \leq \theta^{i}\delta_0 \ \ \Rightarrow \ \ \sup_{x \in B_{\theta^{i+1}}(0)} d(l(x),{_{i+1}l}(x)) \leq \theta^{i+1}(\delta_0+\theta^{-1}d_0).
\]
Thus, the triangle inequality implies (\ref{ssa}) with $\sigma=\theta^{i+1}$ and $c= \delta_0+2\theta^{-1}d_0$ from which  (\ref{ssa}) follows for any $\sigma \in (0,1)$ with a slightly worse constant $c$.

We now sketch how the above analysis implies that if $x_0 \in X$ is an order 1 point of a harmonic map $u:\Omega \to X$, then some neighborhood of $x_0$ maps into $\mathbb{R}$, more precisely into an isometric, totally geodesic embedding of $\mathbb{R}$ in $X$, which we will refer to simply as an apartment.
   It is sufficient to check this property when $u(x_0)$ is the vertex $V$ since every point of $X \setminus V$ has a neighborhood contained in a single apartment.   Since $x_0$ is an order 1 point, a tangent map $l$ of $u$ is a homogeneous degree 1 map effectively contained in an apartment which we denote $X_0$. 

Suppose, for the sake of contradiction, that no neighborhood of $x_0$ maps into $X_0$. Then take $x$ close to $x_0$ and let $B=B_r(x)$ be the largest ball centered at $x$ contained in $u^{-1}(X \setminus X_0)$, and choose $x_1 \in \partial B \cap u^{-1}(X_0)$. By this choice of $x_1$, about half the image of  $B_\sigma(x_1)$ is contained in $X \setminus X_0$  for $\sigma>0$ sufficiently small. However, the same analysis as in the previous paragraph implies (\ref{ssa}) with $l=l_{x_1,V}$ and $x_1=0$. Combined with (\ref{exeffct}) with $l=l_{x_1,V}$,  this implies that $u$ maps most points of $B_{\sigma}(x_1)$ to $X_0$, a contradiction.  Hence, we conclude that every order 1 point  $x_0$ is a regular point of $u$.

Generalizing the above idea for the component map $v$ in the local decomposition $u=(V,v)$ of harmonic map into a DM-complex, we obtain:

\begin{theorem}\label{order1pointsmodelspace}
There exists no order 1 singular point of the component map $v$.
\end{theorem}

This concludes the sketch of the proof of the regularity theorems.

\section{Regularity results for non-locally compact spaces: Teichm\"uller spaces}

%\subsection{Teichm\"uller spaces}\label{technteich}
The local analysis of harmonic maps into the Weil–Petersson completion of Teichm\"uller space presents two main technical challenges. The first  is that the metric near a boundary stratum is not a product, but  only asymptotically a product.  This phenomenon resembles the one for general DM-complexes. 
The second and perhaps the more challenging issue arises from the lack of local compactness and the degenerating geometry of  $\overline{\mathcal T}$.  A new approach  was necessary: neither  Gromov-Schoen nor our previous work on the rigidity of hyperbolic buildings encountered this difficulty, as those settings involved locally compact complexes.

In the following theorems the regular set ${\mathcal R}(u)$ is the set of points in the domain that possess a neighborhood mapping into a single stratum in $\overline{\mathcal T}$.  The singular set ${\mathcal S}(u)$ is the  complement of the regular set.  

%The local analysis of harmonic maps into the Weil–Petersson completion of Teichm\"uller space presents two main technical challenges. The first  is that the metric near a boundary stratum is not a product, but  only asymptotically a product.  This  difficulty is similar to the  one addressed for general DM-complexes. 
%The second and perhaps the more challenging issue arises from the lack of local compactness and the degenerating geometry of  $\overline{\mathcal T}$.  A new approach  was necessary: neither  Gromov-Schoen nor our previous work on the rigidity of hyperbolic buildings encountered this difficulty, as those settings involved locally compact complexes.
%

\begin{theorem}[Regularity theorem I] \label{RegularityTheorem}
Let  $({\mathcal T}, G_{WP})$ denote the Teichm\"uller space of an oriented surface  of genus $g$  and $p$ marked points such that $k=3g-3+p>0$ with the Weil-Petersson metric and  let $(\overline{\mathcal T}, d_{ \overline{\mathcal T}})$ be its metric completion.  If $(\Omega,g)$ is an  $m$-dimensional Lipschitz Riemannian domain and $u:(\Omega,g) \rightarrow( \overline{\mathcal T}, d_{ \overline{\mathcal T}})$ is a harmonic map, then 
\[
\dim_{\mathcal H} \Big({\mathcal S}(u)\Big) \leq m-2.
\]
\end{theorem}

\begin{theorem}[Regularity theorem II] \label{goto0}
Let   $u:(\Omega,g) \rightarrow( \overline{\mathcal T}, d_{ \overline{\mathcal T}})$ be as in Theorem~\ref{RegularityTheorem}.  For any compact subdomain $\Omega_1$ of $\Omega$, there exists a sequence of smooth functions $\psi_i$  with $\psi_i \equiv 0$ in a neighborhood of ${\mathcal S}(u) \cap \overline{\Omega}_1$, $0 \leq \psi_i \leq 1$ and $\psi_i(x) \rightarrow 1$ for all $x \in \Omega_1 \setminus {\mathcal S}(u)$ such that
\[
\lim_{i \rightarrow \infty} \int_{\Omega} |\nabla \nabla u| |\nabla \psi_i| \ d\mu =0.
\]
\end{theorem}

We now briefly describe how to address the technical  issues mentioned before.
%The first technical challenge is addressed through the metric estimates of the Weil-Petersson metric.   
First, the $C^0$ estimates of  \cite{masur}, \cite{daskal-wentworth}, \cite{wolpert}, \cite{yamada} quantify the way in which  the Weil-Petersson metric   is   asymptotically a product  near the boundary of  $\overline{\mathcal T}$. More precisely, let $s=(s_1,\dots, s_n)$ be holomorphic coordinates on the space $\mathcal T_{c}$ parametrizing  nodal surfaces and $t=(t_1, \dots, t_N) \in {\bf C}^N$   denote the plumbing parameters, which analytically describe the collapsing of the curves in $c$.
Define  the real parameters $(r_i, \theta_i)$  by 
\[
r_i=2\pi^2 (-\log|t_i|)^{-\frac{1}{2}},  \ \ \ \ \theta_i=\arg t_i.
\]
Let
\[
{\bf H}_i=\{(r_i,\theta) \in {\bf R}^2:  r_i>0\}, 
\]
denote the ($i$-th copy of) the model space with its Riemannian metric
\[
ds^2= 4dr_i^2+ r_i^6 d\theta_i^2.
\]
This metric is NPC. Its curvature goes to $-\infty$ as $r_i \to 0$.  
Moreover, the metric is incomplete; its completion is obtained by collapsing the axis $r_i=0$ to  a single point.

Let $h$ denote the product metric  on ${\bf H}_1 \times \dots \times {\bf H}_N$ and $G_{wp}$ the Weil-Petersson metric on ${\mathcal T_c}$. Then (cf. \cite{daskal-wentworth}, \cite{wolpert},  and  \cite{yamada}),
\begin{equation} \label{sumio}
G_{WP}- G_{wp}  \oplus  h =  O(|r|^3)  h
\end{equation}
 where 
$|r|:={\left( \sum_{i=1}^N r_i^2\right)}^{1/2}$
is comparable to the  Weil-Petersson distance  to ${\mathcal T_c}$.

To study harmonic maps into Teichm\"uller space, one needs the stronger $C^1$-estimates established in \cite{daskal-meseC1}. These estimates justify formally differentiating the error term in equation~(\ref{sumio}) to obtain asymptotic bounds for the derivatives of $G_{WP}$ in terms of the corresponding derivatives of the product metric. We do not state the estimates here and instead refer the reader to the original source.

%   
%
%
%
%
%%In this paper, we will apply a variation of the  Gromov-Schoen argument with the completion of Teichm\"uller space $\overline{\mathcal T}$ playing the role of a Euclidean building. Since all the degenerating geometry of $\overline{\mathcal T}$ comes from the model space $\overline{\bf H}$,  we will limit our discussion to  $\overline{\bf H}$. 
% 
% 
% The  technical issue above closely resembles the absence of a product structure for the metric near a lower dimensional stratum in a general DM-complex.
 
 We proceed to argue  as for  DM-complexes (cf. (\ref{Vv})): The harmonic map $u:\Omega \rightarrow \overline{\mathcal T}$ can be written near a stratum ${\mathcal T}_c$ as $u=(V, v)$ where $V$ maps into the smooth stratum ${\mathcal T}_c$ and $v$ maps to the normal space $\overline{\bf H}_1 \times \dots \times \overline{\bf H}_n$. 
The proof resembles the one of Theorems~\ref{regtheorem} and~\ref{goto0DM}. 
%As before,  Theorems~\ref{RegularityTheorem} and~\ref{goto0} follow by estimating the codimension of the singular set and proving decay of the harmonic map near the singular set. 
However, there is  the second, more serious issue:
 The non-local compactness and the degenerating geometry of $\overline{\mathcal T}$. 
 
   Since the source of the non-local compactness of $\overline{\mathcal T}$ is solely due to  $\overline{\bf H}$, we will, for the sake of simplicity, concentrate on maps to $\overline{\bf H}$. Recall that the metric completion $\overline{\bf H}$ of the model space ${\bf H}$ with metric $4dr^2+ r^6 d\theta^2$ is obtained by  attaching the boundary line $\{r=0\}$ and identifying it as a single point $P_0$.   Let $u: \Omega \rightarrow \overline{\bf H}$ be a harmonic map.
  The  harmonic map equations  for points not mapping to $P_0$ are
\begin{equation} \label{hmequphi'}
u_r \Delta u_r = 3  u_r^6 |\nabla u_{\theta}|^2 \ \mbox{ and }  \ u_r^4 \Delta u_{\theta} = -6 \nabla u_r \cdot u_r^3 \nabla u_{\theta}.
\end{equation}
Note that the right hand side of each equation  is bounded since $u$ is Lipschitz (cf. Theorem~\ref{2.4.6}). The left hand side of  each equation, though, involves $u_r$.  Thus,  the  equations degenerate as $u_r(x) \rightarrow 0$.  

 The non-local compactness of
  $\overline{\bf H}$ is due to the angular variable $\theta$ that can be arbitrarily large near $P_0$.  This is reflected into the difficulty in defining tangent maps.
Assume that $x_0$ is a singular point, i.e.~$u(x_0)=P_0$. 
For locally compact NPC spaces, the tangent maps typically take values in the tangent cone of the space at the given point. In contrast, for $\overline{\bf H}$, the tangent maps take values in a bouquet of copies of $\overline{\bf H}$ joined at $P_0$. For instance, at an order one point this bouquet is given by $\overline{\bf H} \vee_{P_0} \overline{\bf H}$. Note that the tangent cone of $\overline{\bf H}$ at $P_0$ is isometric to the half-interval $[0,\infty)$ (cf.~\cite{wolpert-Alex}), and the lack of a good notion of an exponential map to relate this space back to the original space $\bf H$ poses a technical difficulty. We will not describe the details here, but note that the issue is resolved  by considering  approximate embeddings of the image of the tangent map  into $\overline{\bf H}$, and by applying the iterative scheme outlined earlier for DM-complexes to establish regularity.
%\\
%\\
%{\bf George, we can provide more details here since  this paper is currently only 17 pages.}
%\\

 This concludes the sketch of the proof of theorems~\ref{RegularityTheorem} and~\ref{goto0}. Recently, Y.~Sun has refined the argument  above to prove that, in fact, higher-order singular points do not exist (cf. \cite{sun}):

\begin{theorem} \label{thm:yitong}
Let $u : \Omega \to \overline{\mathcal T}$ be a harmonic map. If the image of $u$ intersects a stratum $\mathcal T_c$ of $ \overline{\mathcal T}$, then $u$ maps entirely into $\mathcal T_c$.
\end{theorem}

Theorem~\ref{thm:yitong} completes the circle of ideas developed in \cite{daskal-meseHR} by simplifying the original argument.

\section{Infinite energy pluriharmonic maps}\label{sec: infinite energy}
Theorem~\ref{existence} implies existence of a harmonic map assuming that a finite energy map exists.   In practice, though, it may be difficult or even impossible to do so without making additional assumptions. 

In some situations, it is possible to show  existence of an infinite energy harmonic map.
Infinite energy harmonic maps  have appeared, for example,  in the work of Lohkamp and Wolf. Lohkamp \cite{lohkamp} proved the existence of a harmonic map within a given homotopy class of maps between two non-compact manifolds, under a certain simplicity condition. The most important case is when the domain is metrically a product near infinity. Wolf \cite{wolf} studied infinite energy harmonic maps when the domain is a nodal Riemann surface, and applied this to describe degenerations of surfaces in the Riemann moduli space.
Jost and Zuo (cf.~\cite{jost-zuo}) also sketched a proof of the existence of infinite energy pluriharmonic maps from non-compact K\"ahler manifolds to symmetric spaces and buildings.
Mochizuki \cite{mochizuki-memoirs} rigorously constructed pluriharmonic maps into the symmetric space $GL(r,\C)/U(r)$ by proving existence of pluriharmonic metrics on flat vector bundles. These metrics correspond to harmonic maps via the Donaldson-Corlette theorem (cf.\cite{donaldson}, \cite{corlette2}).

 Mochizuki's theory can be extended to other symmetric spaces  and buildings. For example, 
  
\begin{theorem}[Existence of infinite energy pluriharmonic maps \cite{bddm}] \label{theorem:pluriharmonic} 
Let $M$, $X$ and $\rho:\pi_1(M) \rightarrow \mathsf{Isom}(\widetilde X)$ satisfy the following:
\begin{itemize}
\item $M$ is a smooth quasi-projective variety with universal cover $\widetilde M$
\item $X$ is a locally finite irreducible Euclidean building
  \item  $\rho: \pi_1(M) \rightarrow \mathsf{Isom}(X)$ is a homomorphism that does not fix a point at infinity nor a non-empty proper convex subset of $X$.
 \end{itemize}
Then  there exists a unique  $\rho$-equivariant pluriharmonic map $\widetilde u: \widetilde M \rightarrow X$ of logarithmic energy growth. 
 \end{theorem}

Here, {\it logarithmic energy growth} means that the energy of the map grows at most logarithmically.  For precise definition of this notion, we refer to \cite{bddm}. 

The following theorem extends Theorem~\ref{rigiditytheorem1} by removing the assumption of the existence of a finite energy map.

\begin{theorem}Let $M$ be a quasiprojective variety with universal cover $\widetilde M$, $ \Gamma$  the mapping class group of an oriented  surface  $S$ of genus $g$  and  $p$ marked points such that $k=3g-3+p>0$, and let $\rho: \pi_1(M) \rightarrow \Gamma$ a homomorphism that is sufficiently large. Then there exists a unique $\rho$-equivariant  pluriharmonic map $ u: \widetilde M \rightarrow {\mathcal T}$ of logarithmic energy growth (possibly of infinite energy). If, in addition, the (real) rank of $u$ is $\geq  3$ at some point, then  $ u$ is holomorphic or conjugate holomorphic.
\end{theorem}

%\begin{theorem} \label{theorem:teichmuller}
%Let $M$ as above with $\rho:\pi_1(M) \rightarrow \mathsf{Isom}(\widetilde X)$ satisfying the following:
%\begin{itemize}
%\item $\widetilde X$ is the Weil-Petersson completion of a Teichm\"uller space
% \item  $\rho: \pi_1(M) \rightarrow  \mathsf{Isom}(\widetilde X)$ is proper.
% \end{itemize}
%Then  there exists a unique  $\rho$-equivariant pluriharmonic map $\widetilde u: \widetilde M \rightarrow X$ of logarithmic energy growth.   
%\end{theorem}

\section{Concluding remarks}
In this article, we outlined some applications of harmonic maps to rigidity. At this point, a  question arises: what can be said about the rigidity of spaces beyond those with differentiable manifold structure? Can harmonic maps reveal anything new about group actions on more general NPC spaces? Korevaar and Schoen's ambitious program, as outlined in papers \cite{korevaar-schoen}, \cite{korevaar-schoen2}, and \cite{korevaar-schoen3}, sought to initiate such an approach by defining basic analytic concepts such as energy and by proving the existence and Lipschitz regularity of harmonic maps in this generality. However, despite some efforts, there is no significant progress in proving new rigidity theorems via harmonic maps  in this more general context. One exception is \cite[Corollary 1.5.3]{korevaar-schoen2}, which shows that when the domain is a flat torus, harmonic maps to NPC spaces are totally geodesic. The proof relies on an averaging method, which is based on the fact that domain mollification essentially decreases energy. This process is carried out precisely for flat domains, but in general, error terms arising from the domain's metric make the approach more difficult to implement. 

%{\bf Maybe also add something about korevaar-schoen 3, don't remember what they are proving, do you???}

The only other result in this direction is the intriguing theorem of B. Freidin \cite{freidin}, who, building on earlier work of J. Chen \cite{chen}, managed to prove an analogue of the Eells-Sampson formula for maps into general NPC metric spaces and deduce that harmonic maps from manifolds of positive Ricci curvature are constant. Freidin's method involves a very careful analysis of how the Ricci curvature of the domain interacts with the pullback tensor (full energy-momentum tensor) associated to the harmonic map.

Reflecting upon the previous preliminary results, one might conjecture that a further combination of averaging techniques together with careful exploration of the role of symmetries of the spaces involved might yield new rigidity results via harmonic maps. 
 This approach may be even  fruitful in the smooth setting, where conservation laws on spaces with symmetries (such as symmetric spaces or even hyperbolic manifolds) ought to impose restrictions on the harmonic maps between them.

\section*{Acknowledgements}
%\addcontentsline{toc}{section}{Acknowledgements}
The authors would like to thank Rick Schoen and Karen Uhlenbeck for their continued support, Christine Breiner, Ben Dees, Aditya Kumar, for their careful reading and insightful comments, and Miya Mese-Jones for assistance with the figures.


\begin{thebibliography}{99}

\bibitem{abikoff} W. Abikoff. {\it The real analytic theory of Teichm\"uller space.} Springer, Berlin (1980).

\bibitem{almgren}
F.~J.~Almgren, Jr.,
\newblock {\em Almgren's Big Regularity Paper: $Q$-valued functions minimizing Dirichlet's integral and the regularity of area-minimizing rectifiable currents up to codimension two}.
\newblock World Scientific Monograph Series in Mathematics, Vol.~1, World Scientific Publishing Co., Inc., River Edge, NJ, 2000.


%\bibitem{breiner-dees-mese} C.~Breiner, B.~Dees, and C.~Mese.  {\it Harmonic Maps into Euclidean Buildings and Non-Archimedean Superrigidity}. 	arXiv:2408.02783
%
%\bibitem{bader-furman} U.~Bader and A. ~Furman. {\it An extension of Margulis’s superrigidity theorem.} Dynamics. Geometry, Number Theory--The Impact of Margulis on Modern Mathematics, pp. 47-65. University Chicago Press, Chicago, IL (2022).

\bibitem{bddm}
D.~Brotbek, G.~Daskalopoulos, Y. Deng and C.~Mese.
{\it Pluriharmonic maps into buildings and symmetric differentials},
arXiv:2206.11835, 2022.

\bibitem{chen} J. Chen. {\it On energy minimizing mappings between and into singular spaces.} Duke Math. J. 79, (1995) 77-99.

\bibitem{chu}  T. Chuxis, {\it Ph.D Thesis.} Columbia University (1976).

\bibitem{corlette} K.~Corlette. 
{\it Archimedean superrigidity and hyperbolic geometry}.  
Ann. Math. 135 (1990), 165-182.

\bibitem{corlette2}  K.~Corlette. {\it Flat $G$-bundles with canonical metrics}.   J.~Differential Geom. 28 (1988) 361-382.

\bibitem{daskal-meseMRL} G.~Daskalopoulos and C.~Mese.  
{\it Uniqueness of equivariant harmonic maps to symmetric spaces and buildings}.  
Math. Res. Lett. 30 (2023), 1639-1655.

\bibitem{daskal-meseDM} G.~Daskalopoulos and C.~Mese.  
{\it On the singular set of harmonic maps into DM-complexes.}  
Mem. Amer. Math. Soc. 239 (2016).

\bibitem{daskal-meseER} G.~Daskalopoulos and C.~Mese.  
{\it Essential regularity of the model space for the Weil-Petersson metric}.  
J. Pure Appl. Math., published online Aug.~2016.

\bibitem{daskal-meseC1} G.~Daskalopoulos and C.~Mese.  
{\it $C^1$-estimates for the Weil-Petersson metric.}  
Trans. Amer. Math. Soc. 369 (2017), 2917-2950.

\bibitem{daskal-meseHR} G.~Daskalopoulos and C.~Mese.  
{\it Rigidity of Teichm\"{u}ller space.}  Inventiones mathematicae 224 (2001) 1-126.

\bibitem{daskal-meseSR} G.~Daskalopoulos, C.~Mese and A.~Vdovina.  
{\it Superrigidity of hyperbolic buildings.}  
Geom. Funct. Anal. 21 (2011), 1-15.

\bibitem{daskal-wentworth} G. Daskalopoulos and R. Wentworth. {\it Classification of Weil-Petersson isometries} Amer. J. Math. 125 (2003) 941-975. 

\bibitem{donaldson} S.~Donaldson. {\it Twisted harmonic maps and the self-duality equations.} Proc.~London Math.~Soc.~55 (1987) 127-131.

\bibitem{eells-sampson}J.~Eells and J.H.~Sampson, 
\emph{Harmonic mappings of Riemannian manifolds}, 
Amer.~J.~Math.~86 (1964), 109-160.

\bibitem{freidin} B. Freidin. {\it A Bochner formula for harmonic maps into non-positively curved metric spaces.} Calculus of Variations and Partial Differential Equations Vol. 58(121), (2019).

\bibitem{gromov-schoen} M.~Gromov and R.~Schoen.  
\textit{Harmonic maps into singular spaces and $p$-adic superrigidity}.  
Publ. Math. Inst. Hautes Etudes Sci. \textbf{76} (1992), 165-246.

\bibitem{hartman}
P.~Hartman, 
\emph{On homotopic harmonic maps}, 
Canad.~J.~Math.~\textbf{19} (1967), 673-687.
\bibitem{jost-yau} J.~Jost and S.-T.~Yau.  \textit{Harmonic maps and superrigidity. Differential geometry: partial differential equations
on manifolds}  245-280, Proc. Sympos. Pure Math., 54, Part 1, Amer.~Math.~Soc., Providence, RI, 1993.

\bibitem{jost-zuo} J. Jost and K. Zuo.  {\it Harmonic maps of infinite energy and rigidity results for representations of fundamental groups of quasi-projective varieties.} J. of Diff. Geom. 47 (1997) 469-503.

\bibitem{kleiner-leeb}
B.~Kleiner and B.~Leeb, 
\emph{Rigidity of quasi-isometries for symmetric spaces and Euclidean buildings}, 
Publ.~Math.~Inst.~Hautes \'Etudes Sci.~86 (1997), 115-197.


\bibitem{korevaar-schoen} N.~Korevaar and R.~Schoen.  
{\it Sobolev spaces and harmonic maps into metric space targets.}  
Comm. Anal. Geom. 1 (1993), 561-659.

\bibitem{korevaar-schoen2} N.~Korevaar and R.~Schoen.  
{\it Global existence theorem for harmonic maps to non-locally compact spaces.}  
Comm. Anal. Geom. 5 (1997), 333-387.

\bibitem{korevaar-schoen3} N.~Korevaar and R.~Schoen.  
{\it Global existence theorems for harmonic maps: finite rank spaces and an approach to rigidity for smooth actions}.  
Unpublished manuscript.

\bibitem{lohkamp} J.~Lohkamp, \textit{An existence theorem for harmonic maps}. Manuscripta Math 67 (1990) 21-23.

%\bibitem{margulis} G.~Margulis.  
%\emph{Discrete groups of motions on non-positive curvature.}  
%Amer. Math. Soc. Transl. 109 (1977), 35-45.

\bibitem{margulis}
G.A.~Margulis, 
\emph{Discrete subgroups of semisimple Lie groups}, 
Ergebnisse der Mathematik und ihrer Grenzgebiete (3), vol. 17, Springer-Verlag, Berlin, 1991.

\bibitem{masur} H. Masur. {\it Extension of the Weil-Petersson metric to the boundary of Teichm\"uller space.} Duke Math. J., 43(3) (1976) 623-635.

\bibitem{papa} J. McCarthy and A. Papadopoulos. {\it  Dynamics on Thurston's sphere of projective measured foliations}. Comment. Math. Helv.~64 (1989) 133-166.

\bibitem{mok-siu-yeung}
N.~Mok, Y.-T.~Siu, and S.-K.~Yeung,
\emph{Geometric superrigidity},
Invent. Math. 113 (1993) 57-83.

\bibitem{meseUniqueness} C.~Mese.  
{\it Uniqueness theorems for harmonic maps into metric spaces.}  
Commun. Contemp. Math. 4 (2002), 725-750.

\bibitem{mochizuki-memoirs}  T.~Mochizuki.  {\it Asymptotic behaviour of tame harmonic bundles and an application to pure twistor D-modules}. Memoirs of the AMS 185 (2007). 

\bibitem{mostow}
G.D.~Mostow, 
\emph{Strong rigidity of locally symmetric spaces}, 
Annals of Mathematics Studies, No.~78, Princeton University Press, Princeton, N.J., 1973.

\bibitem{reshetnyak} Y.~G.~Reshetnyak.  
{\it Nonexpanding mappings in a space of curvature no greater than $K$.}  
Sib. Mat. Zh. 9 (1968), 918-927.

%\bibitem{schoen-uhlenbeck} R.~Schoen and K.~Uhlenbeck.  
%{\it A regularity theory for harmonic maps.}  
%J. Differential Geom. 17 (1982), 307-335.

%\bibitem{simon} L.~Simon.  
%{\it Regularity Theory for Harmonic Maps.}  
%In: \emph{Theorems on Regularity and Singularity of Energy Minimizing Maps}, Lectures in Mathematics ETH, BirkhÃ€user, 1996.

\bibitem{tromba} A. Tromba. {\it  On a natural algebraic affine connection on the space of almost complex structures and the curvature of Teichm\"uller space with respect to its Weil-Petersson metric.} Manu.~Math. 56 (1986) 475-497.

\bibitem[Sch]{schumacher}  G. Schumacher. {\it The curvature of the Petersson-Weil metric on the moduli space of K\"{a}hler-Einstein manifolds.} In Complex analysis and geometry, Univ. Ser. Math. (1993)  339-354, Plenum, New York.

\bibitem{siu} Y.~T.~Siu.  
{\it The complex analyticity of harmonic maps and the strong rigidity of compact K\"{a}hler manifolds.}  
Ann. Math. 112 (1980), 73-111.

\bibitem{sun} Y.~Sun.  {\it Regularity of harmonic maps into the Teichm\"uller space.} Ph.D.~Thesis, Johns Hopkins University.

\bibitem{wolf} M. Wolf. {\it Infinite energy harmonic maps and degeneration of hyperbolic surfaces in moduli spaces}. J. of Diff. Geom. 33 (1991) 487-539.

\bibitem{wolpert}  S.A. Wolpert.  {\it Geometry of the Weil-Petersson completion of Teichmueller space}. Surveys in Differential Geometry, VIII: Papers in Honor of Calabi, Lawson, Siu and Uhlenbeck, S. T. Yau (editor), International Press (2003).

\bibitem{wolpertPJ} S. Wolpert.  {\it Noncompleteness of the Weil-Petersson metric for Teichm\"{u}ller space}, Pacific J. of Math. 61 (1975) 513-576.

\bibitem{wolpertJDG}  S. Wolpert.  {\it Geodesic length functions and the Nielsen problem.} J. Diff. Geo 25 (1987) 275-296.

\bibitem{wolpert-Alex}  S. Wolpert. {\it Behavior of geodesic-length functions on Teichm\"uller space.} J. of Diff. Geom. 79, (2008), 277-334.

\bibitem{yamada}  S. Yamada. {\it Weil-Petersson Completion of Teichm\"{u}ller spaces.}  Math. Res. Let. 11 (2004) 327-344.
\end{thebibliography}
\end{document}